\documentclass[a4paper,11pt,reqno]{article}
\usepackage{latexsym}
\usepackage[latin1]{inputenc}
\usepackage{amsfonts,amssymb,epsfig,amsmath,amsthm}
\usepackage{setspace}
\usepackage{graphics,subfigure}
\usepackage{makeidx}
\usepackage{multicol}
\usepackage{hyperref}
\usepackage{listings}
\usepackage{enumerate}
\usepackage{color}

\usepackage[a4paper,scale={0.73,0.75},marginratio={1:1},footskip=6mm,headsep=7mm]{geometry}

\numberwithin{equation}{section}

\makeatletter
\newtheorem*{rep@theorem}{\rep@title} \newcommand{\newreptheorem}[2]{%
\newenvironment{rep#1}[1]{%
\def\rep@title{\bf #2 \ref{##1} }%
\begin{rep@theorem} }%
{\end{rep@theorem} } }
\makeatother
\newreptheorem{theorem}{Theorem}
\newreptheorem{lemma}{Lemma}

\def \beq {\begin{eqnarray}}
\def \eeq {\end{eqnarray}}
\def \beqn {\begin{eqnarray*}}
\def \eeqn {\end{eqnarray*}}

\newtheorem{theorem}{Theorem}[section]
\newtheorem{itlemma}[theorem]{Lemma}
\newtheorem{itproposition}[theorem]{Proposition}
\newtheorem{itcorollary}[theorem]{Corollary}
\newtheorem{itremark}[theorem]{Remark}
\newtheorem{itdefinition}[theorem]{Definition}
\newtheorem{itexample}[theorem]{Example}
\newtheorem{itclaim}[theorem]{Claim}
\newtheorem{itfact}[theorem]{Fact}

\newenvironment{fact}{\begin{itfact}\rm}{\end{itfact}}
\newenvironment{claim}{\begin{itclaim}\rm}{\end{itclaim}}
\newenvironment{lemma}{\begin{itlemma}}{\end{itlemma}}
\newenvironment{remark}{\begin{itremark}\rm}{\end{itremark}}
\newenvironment{corollary}{\begin{itcorollary}}{\end{itcorollary}}
\newenvironment{proposition}{\begin{itproposition}}{\end{itproposition}}
\newenvironment{definition}{\begin{itdefinition}\rm}{\end{itdefinition}}
\newenvironment{example}{\begin{itexample}\rm}{\end{itexample}}

\newcommand{\be}[1]{\begin{equation}\label{#1}}
\newcommand{\ee}{\end{equation}}
\newcommand{\bl}[1]{\begin{lemma}\label{#1}}
\newcommand{\br}[1]{\begin{remark}\label{#1}}
\newcommand{\brs}[1]{\begin{remarks}\label{#1}}
\newcommand{\bt}[1]{\begin{theorem}\label{#1}}
\newcommand{\bd}[1]{\begin{definition}\label{#1}}
\newcommand{\bp}[1]{\begin{proposition}\label{#1}}
\newcommand{\bc}[1]{\begin{corollary}\label{#1}}
\newcommand{\bfact}[1]{\begin{fact}\label{#1}.}
\newcommand{\bex}[1]{\begin{example}\label{#1}.}
\newcommand{\ec}{\end{corollary}}
\newcommand{\efact}{\end{fact}}
\newcommand{\eex}{\end{example}}
\newcommand{\el}{\end{lemma}}
\newcommand{\er}{\end{remark}}
\newcommand{\ers}{\end{remarks}}
\newcommand{\et}{\end{theorem}}
\newcommand{\ed}{\end{definition}}
\newcommand{\ep}{\end{proposition}}
\newcommand{\epr}{\end{proof}}
\newcommand{\bpr}{\begin{proof}}
\newcommand{\bcl}[1]{\begin{claim}\label{#1}}
\newcommand{\ecl}{\end{claim}}

\newcommand{\ecs}{\end{corollary}}
\newcommand{\eers}{\end{exercise}}
\newcommand{\eexs}{\end{example}}
\newcommand{\eems}{\end{example}}
\newcommand{\els}{\end{lemma}}
\newcommand{\eles}{\end{lemmaex}}
\newcommand{\ets}{\end{theorem}}
\newcommand{\eds}{\end{definition}}
\newcommand{\eps}{\end{proposition}}
\newcommand{\bi}{\begin{itemize}}
\newcommand{\ei}{\end{itemize}}
\newcommand{\ben}{\begin{enumerate}}
\newcommand{\een}{\end{enumerate}}

\def\vbar{\mathchoice{\vrule height6.3ptdepth-.5ptwidth.8pt\kern-.8pt}
   {\vrule height6.3ptdepth-.5ptwidth.8pt\kern-.8pt}
   {\vrule height4.1ptdepth-.35ptwidth.6pt\kern-.6pt}
   {\vrule height3.1ptdepth-.25ptwidth.5pt\kern-.5pt}}
\def\fudge{\mathchoice{}{}{\mkern.5mu}{\mkern.8mu}}
\def\bbc#1#2{{\rm \mkern#2mu\vbar\mkern-#2mu#1}}
\def\bbb#1{{\rm I\mkern-3.5mu #1}}
\def\bba#1#2{{\rm #1\mkern-#2mu\fudge #1}}
\def\bb#1{{\count4=`#1 \advance\count4by-64 \ifcase\count4\or\bba A{11.5}\or
   \bbb B\or\bbc C{5}\or\bbb D\or\bbb E\or\bbb F \or\bbc G{5}\or\bbb H\or
   \bbb I\or\bbc J{3}\or\bbb K\or\bbb L \or\bbb M\or\bbb N\or\bbc O{5} \or
   \bbb P\or\bbc Q{5}\or\bbb R\or\bbc S{4.2}\or\bba T{10.5}\or\bbc U{5}\or
   \bba V{12}\or\bba W{16.5}\or\bba X{11}\or\bba Y{11.7}\or\bba Z{7.5}\fi}}

\def\B{\cal{B}}

\def \R {{\mathbb R}}

\def \N {{\mathbb N}}
\def \PR {{\mathbb P}}
\def \E {{\mathbb E}}

\def \DD {{\mathbb D}}
\def \s {y}

\def \M {{\cal{M}}}

\newcommand{\Keywords}[1]{\par\noindent 
{\small{\em Keywords\/}: #1}}

\newcommand{\classification}[1]{\par\noindent 
{\small{\em AMS 2010 Classification\/}: #1}}

\newcommand{\ba}[1]{\addtocounter{for}{1} \begin{eqnarray}\label{#1}}
\newcommand{\ea}{\end{eqnarray}}

\def\sqr#1#2{{\vcenter{\vbox{\hrule height .#2pt
                             \hbox{\vrule width .#2pt height#1pt \kern#1pt
                                   \vrule width .#2pt}
                             \hrule height .#2pt}}}}

\def\pmb#1{\setbox0=\hbox{#1}%
   \kern-.025em\copy0\kern-\wd0
   \kern.05em\copy0\kern-\wd0
   \kern-.025em\raise.0433em\box0 }
\def\sqr#1#2{{\vcenter{\vbox{\hrule height.#2pt
     \hbox{\vrule width.#2pt height#1pt \kern#1pt
   \vrule width.#2pt}\hrule height.#2pt}}}}

\def\0{{\bf 0}}

\def\ve{\varepsilon}
\def\s{\sigma}
\def\d{\delta}
\def\l{\lambda}

\def\g{\gamma}

\def\a{\alpha}
\def\th{\theta}
\def\b{\beta}

\def\Th{\Theta}

\def\cal{\mathcal}
\DeclareMathOperator{\Id}{Id}
\DeclareMathOperator{\diag}{diag}

\newenvironment{myenumerate}{%
\begin{list}{\labelenumi}
	{%
	\setlength{\itemsep}{0.4em}%
	\setlength{\topsep}{0.5em}%
	\setlength\leftmargin{2.6em}%
	\setlength\labelwidth{2.15em}%
	\setlength{\labelsep}{0.45em}%
	\usecounter{enumi}%
	}%
	}%
{\end{list}
}

\renewenvironment{enumerate}{
\renewcommand{\theenumi}{\arabic{enumi}}%
\renewcommand{\labelenumi}{{\rm(\theenumi)}}%
\begin{myenumerate}}%
{\end{myenumerate}}

{\end{myenumerate}}

{\end{myenumerate}}

\newenvironment{myitemize}{%
\begin{list}{$\bullet$}% 
 	{%
	\setlength{\itemsep}{0.4em}%
	\setlength{\topsep}{0.5em}%
	\setlength\leftmargin{2.6em}%
	\setlength\labelwidth{2.15em}%
	\setlength{\labelsep}{0.45em}%
	}%
	}%
{\end{list}}

\renewenvironment{itemize}{
\begin{myitemize}}%
{\end{myitemize}}

\def\cT{\mathcal{T}}

\title{Density estimates and short-time asymptotics for a hypoelliptic diffusion process
}

\author{Paolo Pigato    \footnote{%
Department of Economics and Finance, University of Rome Tor Vergata, Via Columbia 2, 00133 Roma, Italy. 
E-mail: \texttt{paolo.pigato@uniroma2.it}
\newline
I am grateful to Vlad Bally, Gianmarco Chinello, Giovanni Conforti, Archil Gulisashvili, Sara Mazzonetto, St\'ephane Menozzi for discussion and suggestions. 
}}

\begin{document}
\date{\today}
\maketitle
\noindent
{\bf Abstract:}
We study a system of $n$ differential equations, each in dimension $d$. Only the first equation is forced by a Brownian motion and the dependence structure is such that, under a local weak H\"ormander condition, the noise propagates to the whole system. We prove upper bounds for the transition density (heat kernel) and its derivatives of any order. Then we give precise short-time asymptotics of the density at a suitable central limit time scale. Both these results account for the different non-diffusive scales of propagation in the various components. Finally, we provide a valuation formula for short-maturity at-the-money Asian basket options under correlated local volatility dynamics.
%We study a system of $n$ differential equations, each in dimension $d$. Only the first equation is forced by a Brownian motion and the dependence structure is such that the noise propagates to the whole system. Assuming a weak H\"ormander condition on the coefficients, we prove upper bounds for the transition density (heat kernel) and its derivatives of any order, Gaussian in the case of bounded diffusion coefficient, log-normal or polynomial in the case of linear-growth diffusion coefficient.  Then we give precise short-time asymptotics of the density of the solution at a suitable central limit time scale.  Both these results account of the different non-diffusive scales of propagation of the solution in the various components. Finally, we provide an application to valuation of short-maturity at-the-money Asian basket options under multi-asset local volatility dynamics and discuss connections with well known results in the literature.

\smallskip

\smallskip
\classification{60H10, 60H30, 60H07, 60J60, 60F05, 91G20
}

\medskip
\begin{center}
\line(1,0){250}
\end{center}

\section{Introduction}

We consider a system of $n$ differential equations, each in dimension $d$. The first one is forced by a Brownian motion; the others do not have a stochastic differential component, but their solution is random due to the dependence structure of the system. We write, for $t\in[0,T]$,
\be{system}
\begin{split}
dX^1_t &= B_1(t,X^1_t,\dots, X^n_t)dt+\s(t,X^1_t,\dots, X^n_t)\circ dW_t\\
dX^2_t &= B_2(t,X^1_t,\dots, X^n_t)dt\\
dX^3_t &= B_3(t,X^2_t,\dots, X^n_t)dt\\
& \vdots \\
dX^n_t &= B_n(t,X^{n-1}_t,X^n_t)dt.
\end{split}
\ee
Here, $W$ is a Brownian motion in $\R^d$ and 
each $X^j_t$ is $\R^d$ valued as well. The initial conditions $X^j_0=\xi_j\in\R^d$ are deterministic for $j=1,\dots, n$.

This system is highly degenerate since $W$ is $d$-dimensional, but $X$ is $nd$-dimensional. The Brownian motion acts only on $X^1$, but under a suitable \emph{H\"ormander condition} on the coefficients the noise propagates in the system from lower to higher coordinates. We recall that the H\"ormander condition is a \emph{hypoellipticity} condition, satisfied if the space where $X$ lives is spanned by the diffusion coefficient and the \emph{commutators} (\emph{Lie brackets}) of diffusion and drift coefficients. 
This hypoellipticity condition ensures existence and regularity of the density \cite{Nualart:06,Shigekawa:04,KusuokaStroock:84}, also referred to as the \emph{heat kernel} of the SDE. 
However, since the noise propagates one-way from lower to higher coordinates, the system displays different time scales in the various components, as shown by Delarue and Menozzi in \cite{DelarueMenozzi:10}.

In this paper we assume a weak H\"ormander condition only at the initial point $\xi$ and that the coefficients are smooth. When they have linear growth, we prove log-normal and polynomial bounds for the density of the solution and for its derivatives of any order. When $\sigma$ is bounded, we prove analogous Gaussian upper bounds.
All of these bounds precisely capture the non-diffusive, multi-scale behavior of the process. 
Then, we obtain short-time \emph{diagonal} asymptotics for the heat kernel, looking at the system at a ``central limit'' time scale, at which a Gaussian limit behavior can be observed.  The covariance of the limit Gaussian depends on the commutators at $\xi$, reflecting the underlying hypoelliptic structure. As an application, we derive a short maturity valuation formula for at-the-money Asian basket options on a correlated multi-asset local volatility model.

 \medskip\noindent
{\bf\textit{Density estimates for hypoelliptic SDEs.}}
General Gaussian upper bounds for the density of a SDE are available under suitable H\"ormander conditions on the coefficients \cite{KusuokaStroock:85}. These are diffusive bounds, meaning that the Gaussian density that bounds the heat kernel from above scales in $t$ as the Brownian motion. This type of estimate gives the exact scale of propagation in the elliptic case, but is not so accurate in the hypoelliptic case. In the direction not spanned by the diffusion coefficient $\s$, but by the commutators, the diffusive estimate does not provide the precise scaling in time of the propagation of the noise.
In a hypoelliptic setting, two-sided Gaussian density bounds in the control-Carath\'eodory metric are given, among others, in \cite{jerison,KusuokaStroock:87}. 
Celebrated works of Ben Arous and L\'eandre \cite{BenArousLeandre:91,BenArousLeandreII:91} provide short-time asymptotics for the heat kernel over the diagonal (see next Section \ref{sec:comparison}, in particular \eqref{exp:dec} for a comparison with our result).
These works suppose a \emph{strong H\"ormander} condition for the coefficients, meaning that the noise propagates in the system through the vector fields of the diffusion coefficient and their commutators, without involving the drift.

If $n\geq 2$, system \eqref{system} cannot satisfy a strong H\"ormander condition, since the $j^{th}$ component of $\sigma$ and its commutators vanish for $j= 2,\dots,n$. On the other hand,
the \emph{weak H\"ormander} condition may hold. In this case, the drift has a key role in the propagation of the noise to all the components. For this reason, the randomness propagates with different speeds in different components, producing a \emph{multi-scale} phenomenon. A similar problem in dimension two was considered by Kolmogorov, who solved explicitly the linear case  in \cite{Kolmogorov}; in possibly non-linear situations, density and tube estimates are provided under weak H\"ormander conditions in  \cite{BallyKohatsu:10,cibelli,pigato2018}. A density estimate for a linear, L\'evy driven system is given in \cite{HuangMenozzi}. Supposing only a weak H\"ormander condition, the support of the density may not be the whole space; bounds for the density in such cases  are given in \cite{cmp}. The small time behavior of a relativistic diffusion and of the circular Langevin diffusion, both satisfying only a weak H\"ormander condition, has been addressed by Franchi in  \cite{franchi,franchi2}. For other examples of
short-time heat kernel asymptotics of weak-H\"ormander-type systems see \cite{BA_BO,BA_PA} and references therein.
Gradient estimates and other functional inequalities  for the heat semigroup of
general Kolmogorov and relativistic diffusions are given in \cite{baudoin2020}.

In the case of smooth coefficients, hypoellipticity ensures existence and smoothness of the density. One can wonder what is the minimal regularity of the coefficients such that a solution to \eqref{system} exists and is unique, and in this case what is the regularity needed in order for the density or a certain derivative of the density to exist. For equation \eqref{system}, this question is highly nontrivial, due to the interplay between the geometry of the system and the regularity of the coefficients. We refer to \cite{chaudruderaynal2017,menozzi3,menozzi2,menozzi1,lpp,Veretennikov} for results in this direction.
Systems with a similar structure and a linear drift have been widely studied, from a more analytical perspective, by Pascucci, Polidoro and coauthors (see \cite{Cinti2008,difrancesco_pascucci_2005,lpp}, the bibliography therein and subsequent work). 
Strong existence and uniqueness of solution for a McKean-Vlasov version of \eqref{system} have recently been proved in \cite{hu2021}.

In \cite{DelarueMenozzi:10}, Delarue and Menozzi consider the ``chained" system \eqref{system}. Under a weak H\"ormander condition uniform on $\R^{nd}$, assuming mild regularity of the coefficients and bounded spectrum for the diffusion matrix, they prove lower and upper Gaussian density bounds, by means of the parametrix method. These bounds account precisely of the non-diffusive time scales of the system. In the present paper, we assume smooth coefficients but only local hypoellipticity. When $\sigma$ is bounded, we obtain a Gaussian upper density bound analogous to the one in \cite{DelarueMenozzi:10}, and the corresponding Gaussian bounds for the derivatives of any order. Moreover, we can relax the boundedness hypothesis on $\sigma$, in this case obtaining (as expected) log-normal or polynomial bounds on the density and its derivatives. Let us also note that under our local hypoellipticity hypothesis, global lower density bounds as in \cite{DelarueMenozzi:10} cannot hold. 

In a hypoelliptic setting, several questions are of interest on the topic of heat kernel asymptotics. 
We have already mentioned the short-time results under strong-H\"ormander conditions in \cite{BenArousLeandre:91,BenArousLeandreII:91} and under weak-H\"ormander conditions in  \cite{franchi,franchi2,BA_BO,BA_PA}. A strictly related applied problem is the asymptotic pricing of Asian-type options, see e.g. \cite{pirjol}. However, a general understanding is still lacking and, to the best of our knowledge, no short-time asymptotics for \eqref{system} are yet known.

% {\color{blue} Such type of expansions have classically been employed to study hypoelliptic diffusions, as we do in the present paper (see e.g. \cite{BenArousLeandre:91,BenArousLeandreII:91}).  More recently they have been used in the development of rough paths theory (see e.g. \cite{Lyons1998}) and in signature based machine learning applications.}
%(see e. g. \cite{ML}).
In order to control the non-degeneracy of the system, we use Malliavin calculus techniques for which we refer to \cite{BC14,Nualart:06}. Our asymptotic results also rely on a stochastic Taylor expansion, that we use to isolate the Gaussian principal term.

 \medskip\noindent
{\bf\textit{Applications.}}
System \eqref{system} models, for example, $n$  coupled oscillators \cite{DelarueMenozzi:10}.  Similar systems have been studied in recent years in fluid dynamics to understand \emph{turbulence}  \cite{MSV07,ABCFP,FGV}. In these models the state variables represent different scales or wave-modes and an energy flow transfers from larger to smaller scales through local interactions. For other examples of application to physical systems, see \cite{DelarueMenozzi:10,Talay} and the references therein.

The same \emph{cascade} structure is used as an auxiliary tool in the study of non-linear Hawkes processes with Erlang memory kernels, motivated by applications in neuroscience \cite{DITLEVSEN2017,Locherbach2019}. Other applications in neuroscience, with $n=2$, are the stochastic Jansen and Rit Neural Mass Model \cite{Ableidinger} and the hypoelliptic stochastic
FitzHugh-Nagumo neuronal model \cite{leonsamson2018,2021Buckwar}.

For $d=1,n=2$, equation  \eqref{system} is used in finance to price Asian options on one asset (cf. \cite{ALZIARY,BPV,pirjol}). In Section \ref{sec:pricing} we consider a correlated multi-asset local volatility model and derive, from our theoretical short-time asymptotics, a valuation formula for short-maturity at-the-money Asian basket options (cf. \cite{DEELSTRA200455,DEELSTRA20102814}). 

 \medskip\noindent
{\bf\textit{Content of the paper.}}
In Section \ref{setting} we present setting and hypotheses and state our results, discuss relations with existing work and applications. In Section \ref{section:global} we prove the upper bound for the density and its derivatives, in Section \ref{section:shorttime}  the short-time asymptotics. 

\section{Setting and results}\label{setting}

In this paper, we assume Lipschitz and linear growth conditions on the coefficients. Therefore, \eqref{system} has a pathwise unique solution, which is a diffusion $X$ in $\R^{nd}$ satisfying
\be{eqn}
X_0=\xi, \quad dX_t=B(t,X_t)dt+\bar{\s}(t,X_t) \circ dW_t,
\ee
with initial condition $\xi=
\left(\begin{smallmatrix}
\xi_1\\
\vdots\\
\xi_n
\end{smallmatrix}\right)
\in\R^{nd}$, where
the coefficients have the form
\be{defsigma}
B= 
\begin{pmatrix}
B_1 \\ \vdots \\ B_n   
\end{pmatrix}
\quad\mbox{ and }\quad
\bar{\s}= 
\begin{pmatrix}
\s \\ \0_d \\ \vdots \\ \0_d   
\end{pmatrix}
\ee
and
\be{depB}
\mbox{for } j>1,\,B_j(t,x)\mbox{ depends only on $t$ and }x_{j-1}, \dots , x_n.
\ee
In \eqref{eqn}, $\circ d W_t$ is the differential in Stratonovic form of a Brownian motion in $\R^{d}$.

For $k,m\in\N$ we denote with $\M(k\times m)$ the set of $k\times m$ matrices with entries in $\R$, and $\M(k)=\M(k\times k)$. Let $M\in \M(k\times m)$. We denote by $\l_*(M)$ (respectively $\l^*(M))$ the smallest (respectively the largest) singular value of $M$. We write $M^T$ for the transposed matrix.
We denote with $\Id_k$ the identity matrix in $\M(k)$, with $\0_k$ the null matrix in $\M(k)$ and with $0_k$ the null vector in $\R^k$. We denote with $\s^i$ and $\bar{\s}^i$ respectively the $i^{th}$ columns of $\s$ and $\bar{\s}$. Note that, for fixed $(t,x)\in \R_+\times \R^{nd}$, we have $\bar{\s}(t,x)\in \M(nd\times d)$; for $i=1,\dots,d$, $\bar{\s}^i(t,x) \in \R^{nd}$; $B(t,x)\in \R^{nd}$. 

Hereafter, $\alpha =(\alpha _{1},...,\alpha _{k})\in \{1,..., nd\}^{k}$ represents a multi-index with length $\left\vert \alpha \right\vert =k$ and $%
\partial_{x}^{\alpha }=\partial _{x_{\alpha _{1}}}...\partial _{x_{\alpha
_{k}}}$. We allow the case $k=0$, giving $\alpha=\emptyset$, and $\partial_x^\alpha f= f$. We write $\nabla$ for the gradient and $\nabla^2 $ for the Hessian matrix of a scalar function. We also write $\nabla_y f(x,y),\,\nabla^2_y f(x,y)$ to denote gradient and Hessian w.r.t. the variable $y$. We denote with $J_{x_j} f(t,x_1,\dots,x_n)$ the Jacobian matrix of a vector field $f$ w.r.t. the $j^{th}$ $d$-dimensional space variable $x_j$. We write $J f$ for the Jacobian matrix w.r.t. the whole space variable $x=(x_1,\dots,x_n)$. We also denote $|z|$ the Euclidean norm of an element $z\in \R^h, \, h\in \N$.

\medskip\noindent
{\bf\textit{Assumptions on the coefficients.}} 
In this paper, we suppose $B,\s\in \cal{C}^\infty$, and we assume:
\bi
\item[($H_{1}$)] Weak H\"ormander condition at $(0,\xi)$:
\[
\l_*\left(J_{x_{n-1}} B_n\dots J_{x_1}B_2 \, \s(0,\xi)\right) =:\l>0.
\] 
\item[($H_{2}$)] There exists a constant $\kappa>0$ such that $ |\xi|+|B(0, \xi)|\leq \kappa$ and for any multi-index $\a$ with $|\a|\geq 1$, for all $y\in \R^{nd}$, for all $t\in [0,T]$,
\[
|\partial_x^\a B(t,y)|+
|\partial_t\partial_x^\a B(t,y)|+\sum_{i=1}^d |\partial_x^\a \s^i(t,y)| +|\partial_t\partial_x^\a \s^i(t,y)|
\leq \kappa.
\]
\ei
We consider also the following (stronger) versions of ($H_{2}$):
\bi 
\item[($H_{2}'$)] Hypothesis ($H_{2}$) holds and for all $y\in \R^{nd}$, for all $t\in [0,T]$,
\[
\sum_{i=1}^d |\s^i(t,y)| \leq \kappa.
\]
\item[($H_{2}''$)] Hypothesis ($H_{2}$) holds and there exist smooth functions 
$\mu_1:\R\times \R^{nd}\to \R^{d}$ and $v:\R\times \R^{nd}\to \M(d)$ 
such that
 the coefficients of the first component can be expressed as
\[
B_1(t,y)=\diag(y_1)\mu_1 (t,y),\quad
\s(t,y)=\diag(y_1) v(t,y),
\]
where $\diag(y_1)\in  \M(d)$ is the diagonal matrix with the entries of the vector $y_1\in \R^d$ on the diagonal, and for all $y\in \R^{nd}$, for all $t\in [0,T]$,
\[
|\mu_1(t,y)|+
\sum_{i=1}^d |v^i(t,y)| \leq \kappa,
\]
where $(v^i)'s$ are the columns of $v$.
\ei
Note that in ($H_{2}'$), the additional hypothesis only concerns the diffusion coefficient and in
($H_{2}''$) it only concerns coefficients of the first component $(X^1)$. 

Let $\lfloor \cdot \rfloor$ denote the integer part function. We set the ``degree'' of $h=1,\dots, nd$ as $g_h=g(h)=2\lfloor \frac{h-1}{d}\rfloor+1$.
We also define the degree of a multi-index $\a$ as
\be{g}
g(\a)=\sum_k g(\a_k).
\ee
\ben
%The 'scale matrix': f
\item[$(D_1)$] For fixed $t>0$, we set
$\cT_t\in \M(nd)$ as a diagonal matrix given by $n$ diagonal blocks in $\M(d)$, with $t^{j-1/2} \Id_d$ as $j^{th}$ diagonal
block, $j=1,\dots, n$. Equivalently, $\cal{T}_t$ can be defined as a diagonal matrix where $(\cal{T}_t)_{h,h}=t^{g_h/2}$ for $h=1,\dots,nd$.
\item[$(D_2)$]
We denote  by $\th=\th(\xi)$ the solution to \eqref{eqn} with a vanishing diffusion coefficient:
\be{eqnd}
\th_0=\xi,\quad d\th_t=B(t,\th_t)dt.
\ee
\een
For any r.v. $G$ in $\R^h$ absolutely continuous w.r.t. the Lebesgue measure, we  write $p_G(y)$ for the density of $G$ at $y$. We use the standard notation $p_t(\xi,y)$ for $p_{X_t}(y)$, the density of $X_t$ at point $y$, with $X$ starting from $\xi$ at time $0$ (the heat kernel of \eqref{eqn}).

\begin{theorem}\label{der:theorem}
 For a fixed $T>0$, let $(X_t)_{t\in(0,T]}$ be the solution to \eqref{eqn} and suppose that the coefficients satisfy $(H_1)$ and $(H_2)$. Then, for any $t>0$, any $y\in \R^{nd}$, $X_t$ admits an infinitely differentiable density $p_t(\xi,y)$. Let $\a$ be a multi-index and $\partial_y^\a p_t(\xi,y)$ the derivative of the heat kernel w.r.t. the second component $y\in \R^{nd}$.

\noindent 
(i) For any $p > 2$ there exist constants $C_{0,p}$  depending on $p,\l,\kappa,T$ and
$C_{\a,p}>0$ depending on $\a,p,\l,\kappa,T$ such that for any $t\in(0,T]$, $y\in \R^{nd}$,
\begin{eqnarray}
\label{densityubpoly}
p_t(\xi,y)
&\leq & \frac{1}{t^{n^2 d/2} }
\frac{ C_{0,p}}{1+ \left| \cal{T}_t^{-1} (y-\th_t) \right|^p },
\\
|\partial_y^\a p_t(\xi,y)|
&\leq &
\frac{1}{t^{(g(\a)+n^2 d)/2}}
\frac{ C_{\alpha,p}}{1+ \left| \cal{T}_t^{-1} (y-\th_t) \right|^p }. \label{derivativeubpoly}
\end{eqnarray}
\noindent 
(ii) If $(H_1)$ and $(H_2')$ hold, there exists $C_0$  depending on $\l,\kappa,T$ and
$C_{\a}$ depending on $\a,\l,\kappa,T$
 such that for any $t\in(0,T]$, $y\in \R^{nd}$,
\begin{eqnarray}
\label{densityub}
p_t(\xi,y)
&\leq & \frac{C_0}{ t^{n^2 d/2}} 
\exp \left( - \frac{\left| \cal{T}_t^{-1} (y-\th_t) \right|^2}{C_0} \right) \\
|\partial_y^\a p_t(\xi,y)|
&\leq &
\frac{C_{\a}}{ t^{(g(\a)+n^2 d)/2}} 
\exp \left( - \frac{\left|\cal{T}_t^{-1}(y-\th_t) \right|^2}{C_{\a}} \right).
\label{derivativeub}
\end{eqnarray}
\noindent 
(iii) If  $(H_1)$ and $(H_2'')$ hold,  there exists $C_0$  depending on $\l,\kappa,T$ and
$C_{\a}$ depending on $\a,\l,\kappa,T$
 such that for any $t\in(0,T]$, $y\in \R^{nd}$,
\begin{eqnarray}
\label{densityublognormal}
p_t(\xi,y)
&\leq & \frac{C_0}{ t^{n^2 d/2}} 
 \exp\left(- \frac{\log^2(  
 \sqrt{t}
 | \cal{T}_t^{-1} (y-\th_t) | 
 )}{C_0 t} \right) \\
|\partial_y^\a p_t(\xi,y)|
&\leq &
\frac{C_{\a}}{ t^{(g(\a)+n^2 d)/2} }
\exp 
\left(- \frac{ \log^2(  
 \sqrt{t}
 | \cal{T}_t^{-1} (y-\th_t) | 
 )}{C_\alpha t} \right).\label{derivativeublognormal}
\end{eqnarray}\end{theorem}
Before stating our second result we need to introduce some more notations.
\ben
%The 'direction matrix': 
\item[$(D_3)$]
We define $A \in \M(nd)$ as a block-diagonal matrix given by $n$ blocks in $\M(d)$, with the matrix product $J_{x_{j-1}} B_j \dots J_{x_1} B_2 \s(0,\xi)$ as $j^{th}$ diagonal block: 
\[
A =
\begin{pmatrix}
\s(0,\xi) & 		\0_d		& 	 \ddots		&   \0_d\\
\0_d			&J_{x_1} B_2 \s(0,\xi) 	& 		\0_d	&   \ddots\\
\ddots 	&		\0_d		&	\ddots	& \0_d   \\        
\0_d	&	\ddots	&	\0_d		&J_{x_{n-1}} B_n \dots J_{x_1} B_2 \s(0,\xi)	\\        
\end{pmatrix}.
\]
This matrix is invertible because of $(H_1)$. 

\item[$(D_4)$] We define $Q$ as a symmetric positive definite block-matrix in $\M(nd)$, given by $n^2$ blocks in $\M(d)$: for $1\leq l,j \leq n$, the block in position $(l,j)$ is 
\[
\frac{\Id_d}{(l+j-1)(l-1)!(j-1)!}.
\]
For $n\in \N\setminus \{0\}$, we also set
\be{qn}
q_n:= (2\pi)^{n/2} 
(\det Q )^{\frac{1}{2d}}=(2\pi)^{n/2}
\frac{\prod_{j=1}^{n-1} j!}{\sqrt{\prod_{j=1}^{2n-1} j!} }.
\ee
\een
\begin{theorem}\label{theorem:epsilon}
Let $(X_t)_{t\in(0,T]}$ be the solution to \eqref{eqn}, and suppose that the coefficients satisfy $(H_1)$ and $(H_2)$. 
Let $p_t(\xi,y)$ be the density of $X_t$ for positive $t$ and $y\in \R^{nd}$. Let $y:[0,T]\rightarrow \R^{nd}$ be a deterministic path  such that
\be{condy}
\lim_{t\downarrow 0} \cT^{-1}_t(y_t-\th_t)=\bar{y}\in\R^{nd}.
\ee
Then,
the following asymptotics hold:
\[
\begin{split}
t^{n^2d/2}  p_t(\xi,y_t) &\xrightarrow{t \downarrow 0} \frac{e^{-\langle(AQA^T)^{-1}\bar{y},\bar{y}\rangle/2}}{q_n ^{d}  |\det A|}, \\
t^{n^2d/2}\cal{T}_t \nabla_y p_t(\xi,y_t)  &\xrightarrow{t \downarrow 0} \frac{e^{-\langle(AQA^T)^{-1}\bar{y},\bar{y}\rangle/2}}{q_n ^{d}  |\det A|}
(AQA^T)^{-1}\bar{y},\\
t^{n^2d/2}  \cal{T}_t \nabla^2_y p_t(\xi,y_t) \cal{T}_t & \xrightarrow{t \downarrow 0} \frac{e^{-\langle(AQA^T)^{-1}\bar{y},\bar{y}\rangle/2}}{q_n ^{d}  |\det A|}
((AQA^T)^{-1} \bar{y} \bar{y}^T (AQA^T)^{-1} -(AQA^T)^{-1} ).
\end{split}
\]
\end{theorem}
 See the discussion around \eqref{commutatorsmatrix} in Section \ref{sec:comparison} for an interpretation of matrices $A$ and $Q$ and the paragraph after
 \eqref{eq:BS:asian} for an interpretation of condition \eqref{condy}.

 \subsection{Comments and comparisons with related work}\label{sec:comparison}
 
The Gaussian density estimate \eqref{densityub} and an analogous lower estimate were obtained in \cite{DelarueMenozzi:10} by means of the parametrix method, assuming uniform hypoellipticity in $\R^{nd}$ and bounded diffusion coefficient. Here, since  $(H_1)$ (weak H\"ormander condition at $(0,\xi)$) does not imply that the support of the solution is the whole $\R^{nd}$, such a global lower estimate cannot be obtained (cf. \cite[Example 2.3]{pigato2018} and \cite{cmp}). 
In the case of unbounded coefficients, we obtain upper estimates \eqref{densityublognormal} and \eqref{densityubpoly}. In \cite{BallyKohatsu:10}, in the two dimensional case, a Gaussian lower bound with the same scales of propagation is proved for possibly unbounded coefficients. One could also consider further relaxations of hypothesis $(H_2)$, for example assuming polynomial growth on the coefficients, provided that some suitable control on the moments of the solution is still available (see \cite{hutz_jentz} for similar applications).

In derivative estimates \eqref{derivativeubpoly}, \eqref{derivativeub}, \eqref{derivativeublognormal}, 
 the scaling pre-factor $t^{-n^2d/2}$ must be multiplied with $t^{-g(\a)/2}$. Recall $g(\a)=\sum_k g(\a_k)$, and $g(h)$ is increasing in $h=1,\dots, nd$. So, when we differentiate w.r.t. higher coordinates, the prefactor in the estimate for the derivative explodes faster as $t\downarrow 0$. Estimates for the first and second order derivatives w.r.t. the first subvector $x_1$ are known for a similar system with linear drift \cite{difrancesco_pascucci_2005}.

Classical Gaussian estimates for the heat kernel \cite{Aronson,jerison,KusuokaStroock:87} are centered in the initial condition $\xi$. Here, as in \cite{DelarueMenozzi:10} and  similarly to \cite{BCP, BallyKohatsu:10, pigato2018}, we need to transport the initial condition via \eqref{eqnd} in $(\th_t)_{t\in[0,T]}$, solution to the unforced equation, which is the center of the random oscillations.
Then the $j^{th}$ coordinate $X^j_t$ propagates around $\th^j_t$ with speed $t^{j-1/2}$. The distance between $\xi$ and $\th_t$ is of order $t$ and in the elliptic setting the random fluctuations have the diffusive scaling $\sqrt{t}$, so the effect of the deterministic transport of the initial condition is negligible. Here, the effect of the deterministic transport is negligible only in $X^1_t$, since in higher coordinates the distance of the $j^{th}$ component $X^j_t$ from $\th^j_t$ at time
$t$ is of order $t^{j-1/2}$. Such non-diffusive scales of propagation are strictly connected with the H\"ormander theorem (cf. Lemma \ref{depBB}) and the one-way propagation of the noise, from lower to higher coordinates
 (see \cite{DelarueMenozzi:10} for a thorough discussion of this phenomenon).
  
 In Theorem \ref{theorem:epsilon} we prove that, when looking at the system near the transported diagonal, at a suitable time scale, the short-time behavior of the heat kernel, its first and second derivatives is Gaussian, with covariance $AQA^T$. Higher order derivatives could be computed with the same method as well. 
 Matrix $Q$ represents the covariance of the high-dimensional Gaussian ``driving" the chain of SDEs in $\R^{nd}$ (cf. \eqref{J} and the \emph{Hilbert matrix}).
Matrix $A$ reflects the geometry of the commutators in the H\"ormander condition ($H_1$):
 the matrix
\begin{equation}\label{commutatorsmatrix}
\left(\bar{\s}, ([\bar{\s}^i,B])_{i=1,\dots,d},\dots, ([\dots [\bar{\s}^i,B],\dots ,B])_{i=1,\dots,d} \right)
\end{equation}
is a block upper triangular matrix, with the same blocks as $A$ on the diagonal (cf. \eqref{abar}). 
Also at this asymptotic level, we see that the heat kernel explodes as $t^{-n^2 d/2}$ close to $\th_t$ (in the sense of \eqref{condy}). The exponent $n^2 d$ depends on the geometric structure (meaning the number of commutators needed to span the whole space) and is the analogous of the integer $Q(x)$ for the diagonal estimates  $p_t(\xi,\xi)$ in the strong-H\"ormander-no-drift case in \cite[Equation (0.4)]{BenArousLeandre:91}. On the other hand, in the same remarkable work \cite{BenArousLeandre:91}, Ben Arous and L\'eandre prove that $p_t(\xi,\xi)$ may display an exponential behavior in $t$. Here, under $(H_2')$ we expect such exponential behavior of $p_t(\xi,\xi)$ from \eqref{densityub}. Indeed, if $B_j(0,\xi)\neq 0_d$ for some $j\in {1,\dots,n}$ (analogous to the requirement in \cite{BenArousLeandre:91} of a non-vanishing drift), taking $y=\xi$ in \eqref{densityub}, we have
\[
\liminf_{t\downarrow 0} \frac{\left| \cal{T}_t^{-1} (\xi-\th_t) \right|}{t^{3/2 - j}} \geq | B_j(0,\xi)| .
\]
If $j\geq 2$, quantity $| \cal{T}_t^{-1} (\xi-\th_t) |$ explodes as $t\downarrow 0$ and we have that
\be{exp:dec}
\limsup_{t\downarrow 0}
t^{2j - 3}  \log   p_t(\xi,\xi) <0 \quad(\mbox{possibly} -\infty),
\ee
with an exponential decay to $0$ of the heat kernel $p_t(\xi,\xi)$.
In \cite[Theorem 0]{BenArousLeandre:91}, a similar exponential decay is shown under the strong H\"ormander condition, if the drift vector field is not contained in the space spanned by the commutators of the diffusion coefficients up to order two. In this case, the diffusive dynamics is too slow to compensate the drift. Similarly, for the dynamics \eqref{system}, under the conditions $B_j(0,\xi)\neq 0_d$ for some $j\in {2,\dots,n}$, the $j^{th}$ random fluctuations are too small to compensate the drift. We also note that for any $j\geq 2$ the decay of $p_t(\xi,\xi)$ in \eqref{exp:dec} is faster than any possible speed of decay in \cite[Theorem 0]{BenArousLeandre:91}.

In our framework, we can only say that the asymptotic diagonal explosion is at least exponential, since
an analogous lower bound cannot possibly hold under the local H\"ormander condition $(H_1)$. Indeed $p_t(\xi,\xi)$ may not even be contained in the support of the density. Consider $\xi=(1,0)$,
\[
X^1_t=1+W_t,\quad X_t^2=\int_0^t (X_s^1)^2 ds.
\]
Clearly $p_t(\xi,\xi)=0$ for any $t>0$. This dynamics fits into the setting of the present paper, once we localize in order to have bounded coefficients (see \cite[Example 2.3]{pigato2018}). In this case, \eqref{exp:dec} holds with $ \log   p_t(\xi,\xi)=-\infty$. However, even in this case, Theorem \ref{theorem:epsilon} states that the heat kernel close to the ``transported diagonal'' $\theta$ converges to a positive constant. Similar considerations apply to the (arithmetic) Asian option under the Black-Scholes model, i.e.
\begin{equation}
\label{eq:BS:asian}
d X^1_t= \sigma X^1_t dW_t+rX^1_t dt,\quad d X_t^2=X_t^1 dt,
\end{equation}
with $\xi=(\xi_1,0)$, $\xi_1 > 0$. Indeed also in this case $ p_t(\xi,\xi) =0$ for any $t>0$, because a.s. $X_t^1>0$ and 
$X_t^2>0$ for any $t>0$. 
In light of these properties of the density, we remark that the asymptotics in  Theorem \ref{theorem:epsilon} hold true close to the transported diagonal $(\theta_t)_t$, since paths $(y_t)_t$ satisfying \eqref{condy} get close to $(\theta_t)_t$ when $t\to 0$. Indeed, when looking at the $h^{th}$ component, \eqref{condy} implies that $(y_t^h-\theta_t^h)/t^{g_h/2}\to \bar{y}_h$. For example, in  the dynamics \eqref{eq:BS:asian}, $\theta_t=(\xi_1 e^{rt}, \xi_1 (e^{rt}-1)/r)$ and a path $(y_t)_t$ satisfies \eqref{condy} if 
\[
\begin{pmatrix}
y_t^1\\
y_t^2
\end{pmatrix}
=
\begin{pmatrix}
\theta^1_t\\
\theta^2_t
\end{pmatrix}
+ 
\begin{pmatrix}
\bar{y}_1 t^{1/2}\\
\bar{y}_2 t^{3/2}
\end{pmatrix}
+ 
\begin{pmatrix}
o( t^{1/2})\\
o( t^{3/2})
\end{pmatrix}
\quad\mbox{ for } t\to 0.
\]
Two sided asymptotic bounds are obtainable if the diffusion is uniformly hypoelliptic, under suitable ``non-vanishing'' assumptions on $B$ close to the initial condition, using the two sided estimates in \cite{DelarueMenozzi:10}. In the case of linear drift and constant volatility, which is Gaussian and computable (as in \cite{Kolmogorov}), one obtains examples of exponential explosion with explicit speed. More precise estimates for $p_t(\xi,\xi)$ in the non-linear setting could be obtained using large deviations techniques. W.r.t. this type of estimates, the asymptotics in Theorem \ref{theorem:epsilon} look at a ``central limit'' regime.  We refer to \cite{gulisashvili,fgp2020} for a description of the difference between such regimes in a different non-diffusive setting and to the bibliography therein for more information about various types of limit regime.

Matrices of commutators similar to \eqref{commutatorsmatrix} are used  in \cite{ BCP, pigato2018} to define non-isotropic norms adapted to the hypoelliptic setting, analogously to $A$ here. In Theorem \ref{der:theorem}, as in \cite{DelarueMenozzi:10}, we do not need the explicit use of the commutators to separate the  different time scales, since every component $X^j$ has its own speed of propagation. Nevertheless, commutators are hidden behind the multi-scale behavior.
The quasi-norm $\|x\|=\sum_{h=1}^{nd} x_h^{1/g(h)}$ (homogeneous w.r.t. to the dilations group) is used in \cite{difrancesco_pascucci_2005,lpp} and the related stream of research to account of the same type of multi-scale phenomenon.

\subsection{Local vol Asian basket options and other modeling applications}\label{sec:pricing}

Asian basket options are financial derivatives on the time-average of a basket of stocks or indices (see \cite[Chapter 8]{DEELSTRA200455}), that we suppose here to be continuously monitored.
Analogously, Asian spread options are derivatives on the difference between time-averages \cite{DEELSTRA20102814}. These derivative products are common in energy markets, where contracts are often priced on the basis of an average price over a certain period. 
As a (formal) corollary of Theorem \ref{theorem:epsilon} we derive now a short maturity valuation formula for at-the-money options of Asian type under local volatility dynamics. To the best of our knowledge Asian basket options have been considered mostly under Black-Scholes models, leading to log-normal distributions, whereas we allow here the richer dynamics generated by a correlated multi-asset local volatility model as e.g. in \cite{BayerLaurence}. Our result generalizes the at-the-money Asian option asymptotics on a single asset under local volatility in \cite[Theorem 6]{pirjol} (in \cite[Theorem 2]{pirjol} out-of-the-money Asian options are also considered, using large deviations theory.)

We consider $d$ stocks following \emph{correlated} local volatility dynamics, with $s_0\in (0,\infty)^d$ initial prices, 
$r\in\R^d$ vector of respective (constant) interest rates (there would be clearly no problem in adding a dividend term) and local volatilities $\Sigma(t,S_t)=(\Sigma_1(t,S_t),\cdots,\Sigma_d(t,S_t))^T$. Then the $d$-dimensional model, with $\diag(\cdot)$ as in $(H_2'')$, is
\[
S_0=s_0, \quad dS_t=\diag(S_t) r dt+\diag(S_t) \diag(\Sigma(t,S_t))  dB_t
\]
where $B=(B_t)_t$ is a correlated $d$-dimensional Brownian motion with correlation matrix $\rho=(\rho_{i,m})_{1\leq i,m\leq d}\in \M(d)$, meaning that $\rho=LL^T$ where $L\in\M(d)$ is a lower triangular matrix with positive diagonal entries (Cholesky decomposition of $\rho$) and $B_t =L W_t$, with $(W_t)_t$ standard $d$-dimensional Brownian motion.

For a fixed vector of weights $w\in \R^d$ we consider the $w$-weighted time-integral average of the $d$ stocks
\[
\alpha_t=\int_0^t w^T S_s ds.
\] 
In case all the components of $w$ are positive, such objects are the key components of Asian basket options, as typically the option is a function of time-average $t^{-1}\alpha_t$, where $t$ is the maturity. If some elements of $w$ are negative (and, typically, assets appearing with opposite signs have strong negative correlation) this represents an Asian spread option.

If $\Sigma_1(0,s_0)>0,\dots \Sigma_d(0,s_0)>0$, the system is locally hypoelliptic (but not globally). Completing $w^T$ with $d-1$ rows to a non degenerate matrix $\bar{w}\in\M(d)$ and using the classic It\^o-Stratonovic conversion, one can apply Theorem \ref{theorem:epsilon} with $n=2$, and
\[
\begin{split}
\sigma(t,y_1)&=\diag(y_1)\diag(\Sigma(t,y_1)) L,\\
B_1(y_1)&=\diag(y_1)r+\frac{1}{2}\sum_{i=1}^d (J_{x_1} \sigma^i)\sigma^i(t,y_1), \quad B_2(y_1)= \bar{w} y_1.
\end{split}
\]
One can see that as $t\to 0$, the density of $\alpha_t$ converges to a Gaussian density, whose mean is determined by \eqref{eqnd} and the coefficients above. We can expand such mean in $t$ as
\[
w^T s_0 t+O(t^2).
\]
The variance of the limit Gaussian density can be deduced from
\[
AQA^T =
\begin{pmatrix}
\s\s^T(0,s_0)
 & 	{\s\s^T(0,s_0) \bar w^T}/{2}		\\
{\bar w \s\s^T(0,s_0)}/{2}	& {\bar w \s\s^T(0,s_0) \bar w^T}/{3}	
\end{pmatrix}
\]
as
\[
\begin{split}
\frac{w^T \s\s^T(s_0)w}{3}
&=
\frac{1}{3}w^T \diag(s_0) \diag(\Sigma(0,s_0))  \rho \diag(\Sigma(0,s_0)) \diag(s_0) w
\\
&= \frac{1}{3}\sum_{1\leq i,m\leq d} \rho_{i,m}
w_i w_m \Sigma_i(0,s_0) \Sigma_m(0,s_0) s_0^i
s_0^m,
\end{split}
\]
expanding by components.
Therefore we deduce from Theorem \ref{theorem:epsilon} that, as $t\to 0$, the density of
\[
\beta_t= t^{-1/2}(t^{-1}\alpha_t-w^T s_0)
\]
converges (along paths as in \eqref{condy}) to such a Gaussian density  and we have (formally) that
\be{}
\frac{1}{t}\int_0^t w^T S_s ds  - w^T s_0
\approx 
N \sqrt{t}
\ee
where $N$ is a centered Gaussian with variance $w^T \s\s^T(s_0) w/{3}$. Setting the strike $K=w^T s_0$ and computing the value of the option on the limit density, denoting $(\cdot)^+$ the positive part, we obtain the following short-maturity approximations for at-the-money Asian basket calls and puts
\[
\begin{split}
\begin{matrix}
Asian\\
basket\,call
\end{matrix}
 &=
E\left[\left(\frac{1}{t}\int_0^t w^T S_s ds-K\right)^+\right]
\approx
 \sqrt{\frac{
t\sum_{1\leq i,m\leq d} \rho_{i,m}
w_i w_m \Sigma_i(0,s_0) \Sigma_m(0,s_0) s_0^i
s_0^m
}{6\pi}}\\
\begin{matrix}
Asian\\
basket\,put
\end{matrix}
 &=
E\left[\left(K-\frac{1}{t}\int_0^t w^T S_s ds\right)^+\right]
\approx
 \sqrt{\frac{t
\sum_{1\leq i,m\leq d} \rho_{i,m}
w_i w_m \Sigma_i(0,s_0) \Sigma_m(0,s_0) s_0^i
s_0^m
}{6\pi}}\end{split}
\]
which reduces  to \cite[Theorem 6]{pirjol} when the basket is composed of only one stock, $\rho=w=1$.

Also note that  if $\Sigma(\cdot)$ is bounded $(H_2'')$ holds with $v(t,y)=\diag(\Sigma(y_1)) L$, and in this case bounds \eqref{densityublognormal} and \eqref{derivativeublognormal}  hold for this multi-asset local volatility model.

With regard to different applications, we also mention that in \cite{MSV07} a special (linear) case of \eqref{system} is shown to converge in long time to a Gaussian stationary measure, with covariance $AQA^T$. Moreover, we note that the short-time Gaussian behavior made explicit in Theorem \ref{theorem:epsilon}, combined with a splitting scheme, may allow to write high-order numerical schemes for the simulation of such hypoelliptic SDEs as in \cite[Section 3]{2020chev_mel_tub} or \cite[Section 3]{2021Buckwar}, motivated both from finance and neuroscience.

\section{Upper bounds for the density and its derivatives}\label{section:global}

We need several notions and results of Malliavin calculus, for which our main reference is \cite{Nualart:06}. 
We use in particular some recent results on density estimates, for which we refer to \cite{BC14}. In Appendix \ref{app-mall} we recall notations and results used in this paper. 

\begin{remark}\label{rem:indexes}
In this section we often deal with block-matrices composed by $d\times d$ blocks and $nd$-dimensional vectors composed by $d$-dimensional sub-vectors, so notations for the indexes are quite involved. We use indexes $i$ and $m$ to denote columns of a matrix $\M(d\times d)$ or $\M(nd,d)$, so we have $i,m\in\{1,2,\dots,d\}$.  
We use indexes $j$ and $l$ to denote sub-vectors in $\R^{d}$ of vectors in $\R^{nd}$, or sub-matrices in $\M(d \times d)$ of matrices in $\M(nd \times d)$, so we have $j,l \in\{1,2,\dots,n\}$. We use the index $h$ to refer to a scalar component of a vector in $\R^{nd}$, so $h\in\{1,\dots,nd\}$. 
\end{remark}

Recall \eqref{eqn} and \eqref{eqnd}.
For fixed $j=1,\dots, n$ we define the following $d(n-j+1)$-dimensional process $X_t^{(j)}$ and its deterministic counterpart $\th_t^{(j)}$:
\be{Y}
X_t^{(j)}=
\begin{pmatrix}
X_t^j\\
\vdots\\
X_t^n
\end{pmatrix}
,\quad\quad
\th_t^{(j)}=
\begin{pmatrix}
\th_t^j\\
\vdots\\
\th_t^n
\end{pmatrix}
.
\ee
Recall \eqref{eqn}, $(D_1)$ and $(D_2)$. We also set, for $t\in (0,T],$
\be{delta} 
\chi_t=\cal{T}_t^{-1}(X_t-\th_t). 
\ee 
In what follows, all the constants in the estimates may depend on the final time horizon $T$ but not on  $t\in(0,T]$.
\bl{distY}
Suppose that the coefficients in \eqref{eqn} satisfy $(H_2)$. Then there exist a constant depending on $k\in \N$;  $p\geq 2$; $j=1,\dots, n$; $T>0$ such that for $t\in(0,T]$
\be{eqmoments}
\|X_t^{(j)}-\th_t^{(j)}\|_{k,p} \leq C_{k,p} t^{j-1/2}.
\ee
Moreover, if $(H_2')$ is also satisfied, there exists $C>0$ such that for $t\in(0,T]$
\begin{equation}\label{eq:tail:per:component}
\PR\Big( 
t^{1/2-j} \sup_{s\in[0,t]}|X_s^{(j)}-\th_s^{(j)}|
\geq a \Big)\leq
 C \exp\left(-\frac{a^2}{C}\right).
\end{equation}
If $(H_2'')$ is satisfied, there exists $C>0$ such that for $t\in(0,T]$
\[
\PR\Big( 
t^{1/2-j} \sup_{s\in[0,t]}|X_s^{(j)}-\th_s^{(j)}|
\geq a \Big)\leq
 C \exp\left(-
\frac{ (\log (a\sqrt{t}))^2}{Ct}
 \right).
\]
\el
\begin{corollary}
Under $(H_2)$, for $k\in \N$, $p\geq 2$, 
\be{momentschi}
\|{\chi_t}\|_{k,p}\leq C_{k,p},
\ee
where $C_{k,p}$ depends on $k,p$ and the constant $\kappa$.  If also $(H_2')$ is satisfied,
\be{tailX}
\PR( |\chi_t |\geq a)\leq
 C \exp\left(-\frac{a^2}{C}\right).
\ee
If instead $(H_2'')$ is satisfied,
\be{tailXlognormal}
\PR( |\chi_t |\geq a)\leq
C \exp\left(- \frac{(\log(a\sqrt{t}))^2}{Ct} \right).
\ee

\end{corollary}
\bpr \emph{(of Lemma \ref{distY})}
Estimate \eqref{eqmoments} for $j=1$ follows from $\|X_t-\xi\|_{k,p} \leq C_{k,p} t^{1/2}$, which is standard under $(H_2)$, and $|\th_t-\xi| \leq C t$.
From \eqref{depB} we have that $B^{(j)}(t,X_t)$ depends only on $t$ and $X^{(j-1)}$. 
Therefore, for  $j\geq 2$, with an abuse of notation,
\be{groit}
X_t^{(j)}-\th_t^{(j)}
=
\int_0^t \big( B^{(j)}(s,X^{(j-1)}_s)-B^{(j)}(s,\th_s^{(j-1)}) \big) ds.
\ee
By induction on $j$, reducing the norm of the time integral to the norm of $X^{(j-1)}_s-\th_s^{(j-1)}$, one gets \eqref{eqmoments}.

Now, noticing
$|X_s-\theta_s|\leq |X_s^{(2)}-\theta_s^{(2)}|+|X_s^1-\theta_s^1|$, we have
\[
|X_t^{(2)}-\theta_s^{(2)}| \leq
\int_0^t | B^{(2)}(s,X_s)-B^{(2)}(s,\theta_s) | ds
\leq
C\Big(\int_0^t |X_s^1-\theta^1_s | ds+ \int_0^t |X_s^{(2)}-\theta^{(2)}_s| ds\Big)
\]
Using Gronwall inequality we get
\[
|X_t^{(2)}-\theta_s^{(2)}| \leq
C
\Big(\int_0^t  |X_s^1-\theta^1_s |  ds\Big)
\]
and as a consequence
\begin{equation}\label{x1x1}
\sup_{s\leq t} |X_s-\theta_s|=\sup_{s\leq t} |X_s^{(1)}-\theta_s^{(1)}| \leq C \sup_{s\leq t} |X_s^{1}-\theta_s^1|
\end{equation}
Under $(H_2')$, one can show (for example as in classic \cite[Proposition 8.1]{bass:book}, plus a Gronwall type argument, recalling also $|\th_t-\xi| \leq C t$) that 
\[
\PR\Big( t^{-1/2} \sup_{s\leq t} |X_s^{1}-\theta_s^1| \geq a\Big)\leq
 C \exp\left(-\frac{a^2}{C}\right).
\]
The same estimate holds with $X_s^{(1)}$ instead of $X_s^{1}$ as a consequence of \eqref{x1x1}. Therefore \eqref{eq:tail:per:component} holds for $j=1$.
 Equation \eqref{groit} implies, with $L$ Lipschitz constant (depending on $\kappa$) in space for $B$,
\[
|X_t^{(j)}-\th_t^{(j)}|
\leq 
L \int_0^t \big|X^{(j-1)}_s-\th_s^{(j-1)}\big| ds \leq 
L t \sup_{s\in [0,t]}\big|X^{(j-1)}_s-\th_s^{(j-1)}\big|
\]
and we can prove \eqref{tailX} by induction on $j$ as before.

Under $(H_2'')$,
\[
X_t^1=\diag(\xi_1) \exp\Big(
\int_0^t \mu_1(s,X_s) ds
+
\sum_{i=1}^d \int_0^t  v^i(s,X_s) \circ dW_s^i
\Big)
\]
and similarly for $\theta^1$. We have
\[
\sup_{s\leq t}|X_s^{1}-\th_s^{1}|
\leq C |\xi|
\exp\Big( \sup_{s\leq t} \Big(Ct +\sum_{i=1}^d \Big|\int_0^t  v^i(s,X_s) dW_s^i \Big| \Big)\Big)
\]
Now the usual estimate using the boundedness of $v$ gives
\[
\PR\Big(  \sup_{s\leq t} |X_s^{1}-\theta_s^1| \geq a \sqrt{t}\Big)\leq
 C \exp\left(-
 (\log (a\sqrt{t}))^2/(Ct)
 \right).
\]
By induction, as before, we prove the last point of the statement.
\epr

For $f,g:\R^{+}\times \R^{m}\rightarrow \R^{m}$ we define the directional derivative (w.r.t. the space variable $x$) $\partial_{g}f(t,x)=\sum_{k=1}^{n} g^{k}(t,x)\partial _{x_{k}}f(t,x)$, and  we recall that the {Lie bracket} (or  {commutator})
is defined as $[g,f](t,x)=\partial _{g}f(t,x)-\partial _{f}g(t,x)$. Let us denote as follows the iterated directional derivative of a vector field $\phi:\R_+\times\R^{nd}\rightarrow \R^{nd}$ w.r.t. $B$:
\be{defitder}
\partial_B^0 \phi =\phi;\quad 
\partial_B^l \phi=\partial_B(\partial_B^{l-1} B)\mbox{ for } l\geq 1. 
\ee
Similarly, we use the following notation for the iterated Lie Brackets: 
\be{defitbra}
\phi= [B,\phi]^{(0)};
\quad
[B,\phi]^{(l)}= [B,[B,\phi]^{(l-1)}] \mbox{ for } l\geq 1.
\ee
We write $\times$ to denote sub-vectors not providing any useful information.
\bl{depBB}
Let $B,\bar{\s}$ as in \eqref{defsigma},\eqref{depB} and $\bar{\s}^i$, $i=1,\dots,d$ the columns of $\bar{\s}$. 
\begin{enumerate}
\item
 For $l=1,\dots,n$, for $l+1\leq j\leq n$, we have that 
\be{depDB}
( (\partial_B)^{l-1} B)_j\mbox{ depends only on } (j-l)^{th}\mbox{ to } n^{th} \mbox{ coordinates.} 
\ee
\item
For any $i=1,\dots, d$,
\be{wholeder}
\partial_{\bar{\s}^i} \partial_{B}^{l-1} B= 
\begin{pmatrix}
 \times \\ 
 \vdots \\
 \times \\
J_{x_l} B_{l+1}  \dots J_{x_1} B_2  \s^i \\ 
 0_d \\ 
 \vdots \\ 
 0_d   
\end{pmatrix}
\begin{array}{c}
  \\ 
 l\,\mbox{ blocks }\\
 \\
 (l+1)^{th}\,\mbox{ block }
\\ 
\\
 n-l-1\,\mbox{ blocks }\\ 
 \\  
\end{array}
\ee
and
 the H\"ormander matrix
\be{abar}
\left(\bar{\s}, ([\bar{\s}^i,B])_{i=1,\dots,d},\dots, ([\dots[\bar{\s}^i,B],\dots ,B])_{i=1,\dots,d} \right)
\ee
is a block upper triangular matrix, with blocks $J_{x_l} B_{l+1}  \dots J_{x_1} B_2 \sigma$ for $l=1,\dots, n$ on the diagonal. So, if $(H_1)$ holds, the weak H\"ormander condition holds at $(0,\xi)$. 
\end{enumerate}
\el
\bpr
Using \eqref{depB}, by induction on $l$ we can prove \eqref{depDB} and that for $l=1,\dots,n-1$ the Jacobian matrix of the iterated directional derivative has the following form
\begin{scriptsize}
\[
\begin{matrix}
\\
   \\ 
l\mbox{ blocks} \\
 \\ 
\\
J(\partial_{B}^{l-1} B)= 
\\
\\
\\ 
\\
\\ 
\\ 
\end{matrix}
\begin{pmatrix}
 \times & \dots & \dots & \dots & \dots & \times\\ 
 \vdots &    & &  & & \vdots \\
\times & \ddots &  &  & & \vdots\\ 
J_{x_l} B_{l+1}  \dots J_{x_1} B_2 & \times & \ddots & \\ 
 0_d & J_{x_{l+1}} B_{l+2}  \dots J_{x_2} B_3 & \times & \ddots\\ 
 \vdots & 0_d &  \quad\ddots\quad  & \times & \ddots\\ 
 0_d   &  \dots & 0_d & J_{x_{n-1}} B_n  \dots J_{x_{n-l}} B_{n-l+1} & \times &\dots 
\end{pmatrix}
\]
\end{scriptsize}
where for $l+1\leq j\leq n$, the blocks in the $j^{th}$ row of the matrix  depend only on $(j-l)^{th}$ to $n^{th}$ coordinates. Recalling the form of $\bar{\s}$ in \eqref{defsigma} and multiplying with $\bar{\s}^i$, we get \eqref{wholeder}. Developing definitions \eqref{defitder} and \eqref{defitbra} we find that
\be{wholebra}
[B,\bar{\s}^i]^{(l)} =(-1)^l
\partial_{\bar{\s}^i} \partial_{B}^{l-1} B+v,\quad 
\mbox{ with }v_{j}=0_d \mbox{ for } j=l+1,\dots,n.
\ee
Now, using the anti-symmetry of commutators, we get \eqref{abar}. 
H\"ormander condition holds because matrix \eqref{abar} is non degenerate under $(H_1)$.
\epr
From the point of view of Malliavin calculus, the importance of the H\"ormander condition lies in the fact that it allows to control the \emph{non-degeneracy} of the Malliavin covariance matrix, thanks to results like the \emph{Norris Lemma} \cite{Norris}. The following is a result of this type, adapted to our specific system. The main improvement w.r.t. \cite[Theorem 2.3.3]{Nualart:06} is that we get a quantitative control on the speed of propagation of the diffusion in the different components, so that after a suitable rescaling the constant in the upper bound does not depend on $t\in[0,T]$. This is a main technical step for the proof of the results in this paper. 
\bl{nondegchain}
For $t>0$, let $\g_{\chi_t}$ be the Malliavin covariance matrix of ${\chi_t}$.
For any $p>2$ there exists $C_p>0$, depending on $\l, \kappa,p,T$, such that
$\E\l_*(\g_{\chi_t})^{-p}<C_p$ for all $t\in(0,T]$.
\el
\bpr
For $M\in \M(k\times m)$ we write $|M|_{Fr}$ the Frobenius norm.
 For two positive semi-definite symmetric matrices $M,\bar{M}$ we write $M \leq \bar{M}$ if $\xi^T M \xi \leq \xi^T \bar{M} \xi$ for all $\xi$. We also recall Remark \ref{rem:indexes}.
Following \cite{Nualart:06}, Section 2.3, we set  $Y_t:=\partial_x X_t$ and $Z_t=Y_t^{-1}$. The following SDEs in Stratonovic form are satisfied:
\be{YZ}
\begin{split}
Y_t&=\Id_{nd} + \sum_{m=1}^d\int_0^t J\bar{\s}^m(s,X_s) Y_s \circ dW^m_s+\int_0^t J B(s,X_s) Y_s ds,\\
Z_t&=\Id_{nd} - \sum_{m=1}^d \int_0^t Z_s J\bar{\s}^m (s,X_s) \circ dW^m_s-\int_0^t Z_s J B(s,X_s)ds.
\end{split}
\ee
Moreover, for $\d<t$, the Malliavin derivative of $X_t$ can be expressed as
\[
D_\d X_t=Y_t Z_\d \bar{\s}(\d,X_\d)=\big(Y_t Z_\d \bar{\s}^m(\d,X_\d)\big)_{m=1,\dots, d}.
\]
We have
\[
D_\d {\chi_t} = D_\d  \cal{T}_t^{-1} (X_t-\th_t) = \cal{T}_t^{-1} Y_t Z_\d \bar{\s}(\d,X_\d).
\]
Recall $(D_1)$. We set now $\tilde{A}\in \M(nd)$ as the block-diagonal matrix 
with blocks $\tilde{A}_{j,j}= (-1)^{j-1} A_{j,j}$ as $j^{th}$ diagonal block.
Remark that $\cT_t$ and $\tilde{A}$ commute. Multiplying by $\Id_{nd} = \cT_t \tilde{A} \,\cT_t^{-1} \tilde{A}^{-1}$, we write
\beq
\g_{\chi_t}=\int_0^t D_\d {\chi_t}(D_\d {\chi_t})^T d\d = \cal{T}_t^{-1} Y_t  \cT_t \tilde{A} \bar{\g}_t \tilde{A}^T \cT_t  Y_t^T \cal{T}_t^{-1},
\eeq
with
\beq\label{gammabar}
\bar{\g}_t=\int_0^t \cT_t^{-1} \tilde{A}^{-1}  Z_\d \bar{\s}(\d,X_\d) \bar{\s}(\d,X_\d)^T Z_\d^T \tilde{A}^{-1,T} \cT_t^{-1} d\d.
\eeq
Remark that
\be{gammabar2}
\g_{\chi_t}^{-1}= 
\cal{T}_t Z_t^T  \cT_t^{-1} \tilde{A}^{-1,T}
\bar{\g}_t^{-1} 
\tilde{A}^{-1} \cT_t^{-1}Z_t \cal{T}_t.
\ee
We have to check the integrability of 
$\l_*(\g_{\chi_t})^{-1}=\l^*(\g_{\chi_t}^{-1})$. Recall that $\l^*(\cdot)$ is a sub-multiplicative norm on the set of matrices, and that for two matrices $M,\bar{M}\in\M(nd)$,  $\l^*(M\bar{M})\leq \l^*(M) \l^*(\bar{M}) $. We have
\be{prodlambda}
\l_*(\g_{\chi_t})^{-1} \leq \l^*( \bar{\g}_t^{-1}) \l^*(\tilde{A}^{-1} \cT_t^{-1} Z_t \cal{T}_t
)^2.
\ee

We now deal with $\l^*( \bar{\g}_\d^{-1})= \l_*(\bar{\g}_\d)^{-1}$. 
For $\phi\in C^2(\R^+\times \R^{nd}, \R^{nd})$, applying Ito's formula \cite[Formula (2.63)]{Nualart:06},
\be{devZphic}
\begin{split}
Z_\d \phi(\d,X_\d)&=
\phi(0,\xi)+\int_0^\d Z_s \sum_{m=1}^d [\bar{\s}^m,\phi](s,X_s) dW_s^m\\
&+
\int_0^\d Z_s \left\{[B,\phi]+\frac{1}{2}\sum_{m=1}^d [\bar{\s}^m, [\bar{\s}^m,\phi]] +\frac{d\phi}{dt}\right\}(s,X_s)ds.
\end{split}
\ee
Here $\frac{d\phi}{dt}$ denotes the derivative with respect to the time variable of $\phi$, and $\xi$ is the initial condition of $X$.
For $f:\R\rightarrow \R^{nd}$ let us denote 
${\bf I}^j_{\d} f(\cdot)=\int_0^\d\dots \int_0^{s_2} f(s_1) ds_1\dots ds_j$ for $j\geq 1$ and ${\bf I}^0_{\d} f(\cdot)=f(\d)$.
We develop $Z_\d\phi(\d,X_\d)$ applying \eqref{devZphic} to $\phi$, then to $[B,\phi]$, then to $[B,\phi]^{(2)}$ and so on until $[B,\phi]^{(n-1)}$. We obtain
\[
\begin{split}
Z_\d \phi(\d,X_\d)&= 
\sum_{j=0}^{n-1} [ B,\phi]^{(j)}(0,\xi)\frac{\d^{j}}{j!}
+\sum_{j=0}^{n-1}{\bf I}^j_{\d} \int_0^{\cdot}  Z_u \sum_{m=1}^d  [\bar{\s}^m,[B,\phi]^{(j)}](u,X_u) dW_u^m\\
&+\frac{1}{2}\sum_{j=0}^{n-1}{\bf I}^j_{\d} \int_0^{\cdot} Z_u \sum_{m=1}^d  [\bar{\s}^m, [\bar{\s}^m,[ B,\phi]^{(j)}]](u,X_u)du\\
&+\sum_{j=0}^{n-1}{\bf I}^j_{\d} \int_0^{\cdot} Z_u \frac{d[ B,\phi]^{(j)}}{d t}(u,X_u)du 
+ {\bf I}_{\d}^{n-1} \int_0^{\cdot} Z_u [B,\phi]^{(n)}(u,X_u)du.
\end{split}
\]
Taking now $\phi=\bar{\s}^i,\,i=1,\dots, d$, we find
\[
Z_\d \bar{\s}^i(\d,X_\d)= 
\sum_{j=0}^{n-1} [ B,\bar{\s}^i]^{(j)}(0,\xi)\frac{\d^{j}}{j!}+ R_\d^i,
\]
with $R^i_\d\in\R^{nd}$ given by
\be{rem}
\begin{split}
R_\d^i &= \sum_{j=0}^{n-1}{\bf I}^j_{\d} \int_0^{\cdot}  Z_u \sum_{m=1}^d  [\bar{\s}^m,[B,\bar{\s}^i]^{(j)}](u,X_u) dW_u^m\\
&+\frac{1}{2}\sum_{j=0}^{n-1}{\bf I}^j_{\d} \int_0^{\cdot} Z_u \sum_{m=1}^d  [\bar{\s}^m, [\bar{\s}^m,[ B,\bar{\s}^i]^{(j)}]](u,X_u)du\\
&+\sum_{j=0}^{n-1}{\bf I}^j_{\d} \int_0^{\cdot} Z_u \frac{d[ B,\bar{\s}^i]^{(j)}}{d t}(u,X_u)du + {\bf I}_{\d}^{n-1} \int_0^{\cdot} Z_u [B,\bar{\s}^i]^{(n)}(u,X_u)du.
\end{split}
\ee
Now, from \eqref{wholebra},
we find
\[
\sum_{j=0}^{n-1} [ B,\bar{\s}^i]^{(j)}(0,\xi)\frac{\d^{j}}{j!}=
\begin{pmatrix}
\s^i(0,\xi)  							&+& V^{1,i}_\d	 	\\
- J_{x_1} B_2 \s^i(0,\xi) \d 	&+& V^{2,i}_\d	\\
J_{x_2} B_3 J_{x_1} B_2 \s^i(0,\xi) \frac{\d^2}{2!} 	&+& V^{3,i}_\d	\\
\vdots 									&+& \vdots	\\        
(-1)^{n-1}J_{x_{n-1}} B_n \dots J_{x_1} B_2 \s^i(0,\xi) \frac{\d^{n-1}}{(n-1)!}  &+& V^{n,i}_\d
\end{pmatrix}
,
\]
where $|V^{l,i}_\d|\leq C\d^l$ for all $l=1,\dots,n$. We write now
\be{eq1}
Z_\d \bar{\s}^i(\d,X_\d)=
\begin{pmatrix}
\s^i(0,\xi)  							&+& V^{1,i}_\d&+& R^{1,i}_\d	 	\\
-J_{x_1} B_2 \s^i(0,\xi) \d 	&+& V^{2,i}_\d&+& R^{2,i}_\d	\\
\vdots 									&+& \vdots &+& \vdots	\\        
(-1)^{n-1} J_{x_{n-1}} B_n \dots J_{x_1} B_2 \s^i(0,\xi) \frac{\d^{n-1}}{(n-1)!}  &+& V^{n,i}_\d &+& R^{n,i}_\d
\end{pmatrix},
\ee
where $R_\d^{l,i}$, $l=1,\dots n$ are the $d$-dimensional sub-vectors of  $R_\d^i$. Now, consider a vector field $\phi(\d,x)=(\phi_1,\dots,\phi_n)(\d,x)$ in $\R^{nd}$, each $\phi_j\in \R^d$. Fix $l\in\{1,\dots,n\}$ and suppose that $\phi_j\equiv 0_d$ for $j\geq l$. Then we have that $[\bar{\s}^i,\phi ]_j=0_d$ for $j\geq l$, and $[B,\phi ]_j=0_d$ for $j\geq l+1$ (from the fact that $J B$ is an upper Hessenberg matrix).  With this in mind,
we apply iteratively \eqref{devZphic} to the terms $Z_u[\cdot,\cdot]$ in \eqref{rem}. 
We obtain that $R_\d^{l,i}$ is of order $\d^{l-1/2}$. More precisely, we have the following bound, for any $i=1,\dots,d$, for any $p\geq 2$: 
\be{intbound1}
\E |R_\d^{l,i}|^p \leq  C_p \d^{p(l-1/2)}.
\ee
To prove this inequality, we use H\"older inequality to estimate the moments of integrals in $d s$ and Burkholder inequality to estimate the moments of integrals in $dW_s$, which appear when we apply \eqref{devZphic}. The precise computations to get from \eqref{rem} to \eqref{intbound1} are long and involved from a notational point of view, but they are quite standard and we leave out the details.

We can write \eqref{eq1} as a matrix product:
\be{devZsigma}
Z_\d \bar{\s}(\d,X_\d)=
\begin{pmatrix}
\s(0,x)  							&+& \bar{R}^1_\d	 	\\
-J_{x_1} B_2 \s(0,\xi) \d 	&+& \bar{R}^2_\d	\\
\vdots 									& & \vdots	\\        
(-1)^{n-1} J_{x_{n-1}} B_n \dots J_{x_1} B_2 \s(0,\xi) \frac{\d^{n-1}}{(n-1)!}  &+& \bar{R}^n_\d
\end{pmatrix}.
\ee
with $\bar{R}^{l}_\d=[V^{l,1}_\d+R^{l,1}_\d, \dots, V^{l,d}_\d +R^{l,d}_\d]$. Note that, for fixed $i=1,\dots,d$, $R^i_\d$ in \eqref{rem} is a vector in $\R^{nd}$, whereas here the reminder $\bar{R}^l_\d$, for $l=1,\dots, n$, is a matrix in $\M(d)$. From \eqref{intbound1}, for all $p\geq 2$
\be{intbound2}
\E |\bar{R}_\d^{l}|_{Fr}^p \leq  C_p \d^{p(l-1/2)}
\ee
(the norm used here is the Frobenius norm). Recall that $\tilde{A}$  is non degenerate because of $(H_1)$. 
Using the block-diagonal structure of $\tilde{A}$ we get
\be{AZ}
\cT_t^{-1} \tilde{A}^{-1} Z_\d \bar{\s}(\d,X_\d)= \cal{T}_t^{-1} 
\begin{pmatrix}
\Id_d\\
\Id_d\, \d\\
\vdots\\
\frac{\Id_d}{(n-1)!}\d^{n-1}    
\end{pmatrix}
+
\cT_t^{-1} \tilde{A}^{-1} 
\begin{pmatrix}
\bar{R}^1_\d	 	\\ \vdots	\\        
\bar{R}^n_\d
\end{pmatrix}= t^{-1/2} 
\begin{pmatrix}
\Id_d\\
\Id_d\, \frac{\d}{t}\\
\vdots\\
\frac{\Id_d}{(n-1)!}\left(\frac{\d}{t}\right)^{n-1}    
\end{pmatrix}
+
\tilde{R}_\d
\ee
where 
\[
\tilde{R}_\d
=
\begin{pmatrix}
\tilde{R}^1_\d\\
\vdots\\
\tilde{R}^n_\d    
\end{pmatrix},\quad\quad \tilde{R}^l_\d= t^{-l+1/2}\tilde{A}_{l,l}^{-1} \bar{R}^l_\d,\quad l=1,\dots, n.
\]
From \eqref{intbound2} follows that for any $p\geq 2$: 
\be{intbound}
\E \left| \int_0^{\d} \tilde{R}_s^l (\tilde{R}^j_s)^T
ds \right|^p_{Fr} \leq  C_p \frac{\d^{p(l+j)}}{t^{p(l+j-1)}}
\ee
For fixed $\ve$ let $\rho$ be the following time-dependent matrix:
\[
\rho^\ve_\d=(\rho^\ve_{l,j}(\d))_{1\leq l,j \leq n}
=
\left(
\left(\frac{1}{\ve} \right)^{l+j-1}
\int_0^\d \tilde{R}_s^l (\tilde{R}^j_s)^T 
ds \right)_{1\leq l,j\leq n}
\]
From \eqref{intbound}, for any $p\geq 2$ we have 
\be{bl}
\E\l^*(\rho^\ve_{t \ve})^p\leq C\sup_{l,j}  \frac{(t \ve)^{p(l+j)}}{(t\ve)^{p(l+j-1)}}\leq C (t \ve)^p.
\ee

Recall now $(D_4)$. We have $\det Q = \left((\prod_{j=1}^{n-1} j!)^2 \big/ (\prod_{j=1}^{2n-1} j!)\right)^d $. We introduce the stopping time
\be{defSc}
S_\ve = \inf  \left\{ s \geq 0: \l^*(\rho^\ve_s) \geq \frac{{\l_*(Q)}}{4} 
 \right\} \wedge t,
\ee
For $p\geq 2$ we have, from Markov inequality, the fact that { $\l^*(\rho^\ve_t)$ is increasing in $t$} and \eqref{bl},
\be{estSc}
\PR(S_\ve < t\ve)
\leq
\PR\left(
\l^*(\rho^\ve_{t \ve}) \geq \frac{1}{4} {\l_*(Q)} 
\right)
\leq
\frac{4^p \E \l^*(\rho^\ve_{t \ve})^p}{ {\l_*(Q)}^p }
\leq
C (t \ve)^{p} 
\leq \ve^{p-1}
\ee
for $t\leq T,\,\ve \leq \ve_p$ depending on $p$ and $\kappa$. 
We work on $\big\{S_\ve\geq \ve t \big\}$. We recall that inequality 
\[
(M+\bar{M})(M+\bar{M})^T \geq \frac{1}{2} M M^T- \bar{M} \bar{M}^T
\]
holds for any matrix $M,\,\bar{M}$.
We apply\eqref{AZ}, and obtain
\[
\begin{split}
\bar{\g}_t 
&=\int_0^t \cT_t^{-1} \tilde{A}^{-1} Z_\d \bar{\s}(\d,X_\d) \bar{\s}(\d,X_\d)^T Z_\d^T \tilde{A}^{-1,T}\cT_t^{-1}  d\d\\
&\geq \int_0^{S_\ve} \cT_t^{-1} \tilde{A}^{-1} Z_\d \bar{\s}(\d,X_\d) \bar{\s}(\d,X_\d)^T Z_\d^T \tilde{A}^{-1,T}\cT_t^{-1} d\d\\
&\geq \int_0^{S_\ve} \frac{1}{2} t^{-1} \left( \begin{array}{c}
\Id_d\\
\Id_d\, \frac{\d}{t}\\
\vdots\\
\frac{\Id_d}{(n-1)!}\left(\frac{\d}{t}\right)^{n-1}    
\end{array}\right)
\left( \begin{array}{c}
\Id_d\\
\Id_d\, \frac{\d}{t}\\
\vdots\\
\frac{\Id_d}{(n-1)!}\left(\frac{\d}{t}\right)^{n-1}    
\end{array}\right)^T d\d - \int_0^{S_\ve} \tilde{R}_\d \tilde{R}^T_\d
d\d.
\end{split}
\]
Direct computations give
\[
\begin{split}
\bar{\g}_t 
&\geq \left( \frac{1}{2}Q_{l,j}
\left(\frac{S_\ve}{t} \right)^{l+j-1} 
- \int_0^{S_\ve} \tilde{R}_\d^l (\tilde{R}^j_\d)^T
d\d \right)_{1\leq l,j\leq n}\\
&= \left(\left( \frac{1}{2}Q_{l,j}
- \left(\frac{t}{S_\ve} \right)^{l+j-1} 
\int_0^{S_\ve} \tilde{R}_\d^l (\tilde{R}^j_\d)^T
d\d \right) \left(\frac{S_\ve}{t} \right)^{l+j-1} \right)_{1\leq l,j\leq n}.
\end{split}
\]
Since $\ve\leq \frac{S_\ve}{t}\leq T$, using $\l^*(M)^{-1}=\l_*(M^{-1})$ and $\l^*(M\bar{M})\leq \l^*(M)\l^*(\bar{M})$ we obtain
\[
\begin{split}
\l_*(\bar{\g}_t)&\geq \left[ \frac{1}{2} \l_*(Q)
- \l^* \left(\left(\left(\frac{t}{S_\ve} \right)^{l+j-1} 
\int_0^{S_\ve} \tilde{R}_\d^l (\tilde{R}^j_\d)^T
d\d \right)_{1\leq l,j\leq n}\right)\right]\left(\frac{S_\ve}{t} \right)^{2n-1}
\end{split}
\]
and
\[
\l^* \left(\left(\left(\frac{t}{S_\ve} \right)^{l+j-1} 
\int_0^{S_\ve} \tilde{R}_\d^l (\tilde{R}^j_\d)^T
d\d \right)_{1\leq l,j\leq n}\right)
\leq
\l^*\left(\left(\left(\frac{1}{\ve} \right)^{l+j-1} 
\int_0^{S_\ve} \tilde{R}_\d^l (\tilde{R}^j_\d)^T
d\d \right)_{1\leq l,j\leq n}\right).
\]
Using \eqref{defSc}, since we are on $\{S_\ve\geq t\ve\}$,
\[
\l_*(\bar{\g}_t)\geq \frac{\l_*(Q)}{4} \left(\frac{S_\ve}{t} \right)^{2n-1} \geq 
\frac{\l_*(Q)}{4} \ve^{2n-1}.
\]
Recal \eqref{estSc}. We have that for any $p$ exists $C_p,\ve_p$ such that for any $\ve\leq\ve_p$, $t\leq T$, $|\xi|=1$
\[
\PR(\langle \bar{\g}_t\xi, \xi\rangle< \ve^{2n}) \leq
\PR[S_\ve < t\ve]
 \leq C_p\ve^{p-1}.
\]
The following Lemma is a slight modification of \cite[Lemma 2.3.1]{Nualart:06}. 

\bl{lemma:matrixmoment} 
Let $\g\in \M(n)$ be symmetric, nonnegative definite. We assume that for fixed $p\geq 2$, $\E|\g|^{p+1}_{Fr} < \infty$, and that exists $\ve_0>0$ such that for $\ve \leq \ve_0$,
\[
\sup_{|\xi|=1} \PR [\langle\g \xi,\xi \rangle<\ve]\leq \ve^{p+2n}
\]
Then there exist a constant $C$ such that
\[
\E \l_*(\g)^{-p} \leq C\E|\g|_{Fr}^{p+1} \ve^{-p}_0.
\]
\el
 This implies that for any $p\geq 2$ there exists $C>0$ such that  
 \be{redest}
\E\l_*(\bar{\g}_t)^{-p}\leq C.
\ee 
Now we look at 
$\l^*(\tilde{A}^{-1} \cT_t^{-1} Z_t \cal{T}_t)$.
The matrix $\tilde{A}$ is non-degenerate and does not depend of $t$, so we can as well consider only $\l^*(\cal{T}_t^{-1} Z_t \cal{T}_t
)$. Recall \eqref{YZ}:
\[
Z_t=\Id_{nd} - \sum_{m=1}^d \int_0^t Z_s J\bar{\s}^m (s,X_s) \circ dW^m_s-\int_0^t Z_s J B(s,X_s)ds.
\]
We have
\[
\int_0^t Z_s J \bar{\s}^m(s,X_s) \circ dW^m_s 
 = 
(\mu_{l,j}^m)_{1\leq l,j\leq n}
\]
with $\E|\mu_{l,j}^m|^p_{Fr}\leq Ct^{p(l-1/2)}$. To prove this, we apply \eqref{devZphic} taking for $\phi$ the columns of $J \bar{\s}^m$, for fixed $m=1,\dots ,d$. Then, we apply again \eqref{devZphic}
 to the new terms $Z_u[\cdot,\cdot]$, and iterate the procedure.
Taking for $\phi$ the columns of $J B$ we find
\[
\int_0^t Z_s J B(s,X_s) ds 
 = (\mu_{l,j}^B)_{1\leq l,j\leq n}
\]
with $\E|\mu_{l,j}^B|^p_{Fr}\leq Ct^{p((l-j)\vee 1)}$. Now recall $(D_1)$: $\cal{T}_t$ is a diagonal block-matrix, and the $j^{th}$ block is $(\cal{T}_t)_{j,j}=\Id_d t^{j-1/2}$. Therefore for any random matrix $\mu=(\mu_{l,j})_{1\leq l,j\leq n}$,
\[
\E|(\cal{T}_t^{-1}\mu\cal{T}_t)_{l,j}|^p_{Fr} =t^{p(-l+1/2)}
\E|\mu_{l,j}|^p_{Fr}
t^{p(j-1/2)} =
t^{p(j-l)} 
\E|\mu_{l,j}|^p_{Fr}.
\]
Therefore
\[
\E\l^*( \cal{T}_t^{-1}(\mu_{l,j}^m)_{1\leq l,j\leq n}\cal{T}_t )^p \leq Ct^{p/2}
\quad\mbox{and}\quad\E\l^*( \cal{T}_t^{-1}(\mu_{l,j}^B)_{1\leq l,j\leq n} \cal{T}_t )^p \leq C.
\]
This implies that for any $p$ exists $C$ such that $\E \l^*(\cT_t^{-1} Z_t \cal{T}_t
)^q\leq C$. The proof is concluded once we recall
\eqref{prodlambda}
and
\eqref{redest}.
\epr
We have now all the tools we need to prove Theorem \ref{der:theorem}. 
\bpr \emph{(Theorem \ref{der:theorem})}
Recall $(D_1)$, $(D_2)$, \eqref{g}, \eqref{delta}.
The H\"ormander condition $(H_1)$ and the assumption of infinitely differentiable coefficients imply that $X_t$ is absolutely continuous on $\R^{nd}$ and the density infinitely differentiable, following \cite{Nualart:06}. This can be also seen as a direct consequence of Lemma \ref{deri} and Lemma \ref{nondegchain}. We change variable via $\cT_t z+\th_t=y$.
Since $\cal{T}_t$ and $\th_t$ are deterministic, ${\chi_t}$ is absolutely continuous as well and
\[
p_t(\xi,y)=p_{X_t}(y)=p_{\chi_t} (\cal{T}_t^{-1} (y-\th_t)) \frac{1}{|\det \cal{T}_t|}
=t^{-n^2 d/2} p_{\chi_t} \left( \left( \frac{(y-\th_t)_h}{t^{g_h/2}}\right)_{h=1,\dots, nd} \right) .
\]
Moreover, $p_{\chi_t}(z)$ is infinitely differentiable. 
Applying the chain rule we obtain, for any multi-index $\a\in \{0,\dots,nd\}^k$,
\[
\partial_y^\a p_t(\xi,y)
=t^{-(g(\a)+n^2 d)/2} \partial_z^\a p_{\chi_t} \left( \cal{T}_t^{-1}(y-\th_t) \right).
\]
Using Markov inequality we get the tail estimates from the moment bound \eqref{momentschi}: for any $p\geq 2$, there exists $C_p$ such that
\[
\PR( |\chi_t |\geq z)\leq
 \frac{C_p}{1+|z|^p}.
\]
We apply now Lemma \ref{deri} with $F={\chi_t}$. This tail estimate, together with Lemma \ref{nondegchain}
and
\eqref{momentschi}, implies that for any multi-index $\alpha$ and $p > 2$ there exists $C_{\a,p}>0$ such that
\[
|\partial_z^\a p_{\chi_t} (z)|\leq
  \frac{C_{\a,p}}{1+|z|^p},
\]
and \eqref{densityubpoly} and  \eqref{derivativeubpoly} follows. If now we suppose $(H_2')$ we have the Gaussian tail estimate  \eqref{tailX} and \eqref{densityub} and  \eqref{derivativeub}  follow in the same fashion.
Analogously,  if we suppose $(H_2'')$ we have the log-normal tail estimate  \eqref{tailXlognormal} and \eqref{densityublognormal} and  \eqref{derivativeublognormal}  follow.
\epr

\section{Short-time asymptotics}\label{section:shorttime}

We write now the \emph{Stochastic Taylor development} of $X$, which is the key step to obtain the estimate in Theorem \ref{theorem:epsilon}. Let us introduce the following\\
\noindent {\bf Condition (R)}. Let $\beta=(\beta_t)_{t\in[0,T]}$ be a process in $\R^{nd}$ with
\[
\beta_t=
\begin{pmatrix}
\beta_t^{1}\\
\vdots\\
\b_t^{n} 
\end{pmatrix}
\]
We say that $\b$ satisfies condition {\bf (R)} if for any $k\in\N$, $p\geq 2$, $j=1,\dots ,n$ there exists a constant $C_{k,p}$ depending on $T,\kappa$ in $(H_2)$,
 such that
\be{RR}
\|\b_t^j\|_{k,p}\leq C_{k,p} t^j,\quad\forall t\in(0,T].
\ee
\bl{decomposition}
Let $X$, $\th$, $A$ be given in \eqref{eqn}, $(D_2)$, $(D_3)$. For $ t\in(0,T]$ we have the following decomposition:
\[
X_t-\th_t=A N_t+R_t ,
\]
where $N_t$ is a r.v. in $\R^{nd}$ defined as follows: for $j=1,\dots,n$, $N_t^j$ is the following r.v. in $\R^d$:
\be{J}
N_t^j=
\int_0^t \frac{(t-s)^{j-1}}{(j-1)!} dW_s.
\ee
Moreover, $R_t$ satisfies condition {\bf (R)}.
\el

\bpr
In the following proof, we denote with $\cal{V}^{(m)}$ random vectors in $\R^{nd}$; the superscript $(m)$ is used to denote different vectors, each in $\R^{nd}$. We write
$(\cal{V}^{(m)})^j$ to refer to the $j^{th}$ sub-vector in $\R^d$ of the vector $\cal{V}^{(m)}$.

For $m=1,\dots, n$ we define $\bar{N}_t^{(m)}$ as a r.v. in  $\R^{nd}$ via:
\[
\begin{split}
(\bar{N}_t^{(m)})^j &=N_t^j   \quad\,\mbox{ for } j=1,\dots, m,\\
&=0 \quad\quad  \mbox{ for }j=m+1,\dots, n. 
\end{split}
\]
We also set
\[
\B_t^{(m)}=
\left( 
\begin{array}{c}
(\B_t^{(m)})^{1}\\
\vdots\\
(\B_t^{(m)})^{n} \\
\end{array}
\right)
\]
with
\[
(\B_t^{(m)})^{j}=\int_0^t \dots \int_0^{s_{m\wedge j -1}}  \big(
(\partial_B^{m\wedge j-1} B)_j(s_{m\wedge j} ,X_{s_{m\wedge j}})
-(\partial_B^{m\wedge j-1} B)_j (s_{m\wedge j},\th_{s_{m\wedge j}})
\big)
 d s_{m\wedge j} \dots ds_1
\]
We prove by induction that for any $m=1,\dots, n$
\be{indB}
X_t-\th_t= A \bar{N}_t^{(m)}+\B^{(m)}_t+\bar{ R}_t^{(m)},
\ee
with $(\bar{ R}_t^{(m)})_{t\in (0,T]}$ satisfying condition {\bf (R)}.
Writing $X_t-\th_t$ using \eqref{eqn} and \eqref{eqnd} we obtain
\[
\begin{split}
X_t-\th_t&= \int_0^t  \bar{\s}(s_1,X_{s_1}) \circ dW_{s_1} + \int_0^t  \big(B(s_1,X_{s_1})-B(s_1,\th_{s_1}) \big) ds_1\\
&= \bar{\s}(0,\xi) W_t + \int_0^t \big( B(s_1,X_{s_1})-B(s_1,\th_{s_1})\big) ds_1 +
\int_0^t \big( \bar{\s}(s_1,X_{s_1})-\bar{\s}(0,\xi)\big) \circ dW_{s_1}.
\end{split}
\]
From $(H_2)$, a standard estimate gives that the stochastic integral satisfies condition {\bf (R)}. Moreover 
$\B^{(1)}_t=\int_0^t \big( B(s_1,X_{s_1})-B(s_1,\th_{s_1}) \big) ds_1$, 
so \eqref{indB} is proved for $m=1$. Now we suppose that it holds for $m$ and prove it for $m+1$. Remark that for $j\leq m$,
$(\B_t^{(m)})^{j}=(\B_t^{(m+1)})^{j}$. For $j\geq m+1$,
\[
(\B_t^{(m)})^{j}=\int_0^t\dots \int_0^{s_{m -1}}  
\big(
(\partial_B^{m-1} B)_j(s_{m} ,X_{s_{m}})
-(\partial_B^{m-1} B)_j(s_{m},\th_{s_{m}}) 
\big)
d s_{m} \dots ds_1.
\]
We develop this term using \eqref{eqn}:
\[
\begin{split}
(\B_t^{(m)})^{j} &=
\sum_{i=1}^d 
\int_0^t \dots \int_0^{s_{m}}  
(\partial_{\bar{\s}^i} \partial_B^{m-1} B)_j(s_{m+1} ,X_{s_{m+1}})
\circ dW^i_{s_{m+1}} d s_{m} \dots ds_1
\\
&+ \int_0^t \dots \int_0^{s_{m}}  
\Big(
(\partial_B^{m} B)_j(s_{m+1} ,X_{s_{m+1}})-(\partial_B^{m} B)_j(s_{m+1} ,\th_{s_{m+1}})
\Big)
ds_{m+1} d s_{m} \dots ds_1
\\
&+ \int_0^t \dots \int_0^{s_{m}}  
\Big( (\partial_t \partial_B^{m-1} B)_j(s_{m+1} ,X_{s_{m+1}}) - (\partial_t \partial_B^{m-1} B)_j(s_{m+1} ,\th_{s_{m+1}})
\Big)
ds_{m+1} d s_{m} \dots ds_1\\
&=(U^{(m)}_t)^{j}+(V^{(m)}_t)^{j}+(Z^{(m)}_t)^{j}
\end{split}
\]
We also set $(U^{(m)}_t)^{j}=(V^{(m)}_t)^{j}=(Z^{(m)}_t)^{j}=0_d$ for $j\leq m$, and consequently $U^{(m)}_t,V^{(m)}_t,Z^{(m)}_t$ random vectors in $\R^{nd}$.

From \eqref{wholeder}, $(U^{(m)}_t)^{j}=0$ for $j\geq m+2$. For $j=m+1$ we have 
\[
(U^{(m)}_t)^{m+1}=
\int_0^t \int_0^{s_1}\dots \int_0^{s_{m}}  
J_{x_{m}} B_{m+1} \dots J_{x_1} B_2 \s
(s_{m+1} ,X_{s_{m+1}})
\circ dW_{s_{m+1}} d s_{m} \dots ds_1.
\]
Freezing the integrand in $(0,\xi)$, noting that for $m\geq 1$, using Stochastic Fubini,
\[
\int_0^t \int_0^{s_1}\dots \int_0^{s_{m}}  
 dW_{s_{m+1}} d s_{m} \dots ds_1=
N^{m+1}_t 
\]
and recalling $(D_3)$, the definition of $A$,
we find
\be{Adev}
\| (U^{(m)})^{m+1}- A_{(m+1,m+1)} N^{m+1}_t \|_{k,p}\leq C_{k,p} t^{m+1}
\ee
for all $k,p$. 
This is a classical estimate for the Sobolev norm, and follows from $(H_2)$, H\"older and Burkholder inequalities. Therefore we can write
\[
A \bar{N}^{(m)}_t+U^{(m)}_t=
A \bar{N}^{(m+1)}_t+L^{(m)}_t,
\]
with $L^{(m)}_t$ satisfying condition {\bf (R)}. 
Concerning the second summand, for $j\geq m+1$
\be{Bdev}
(V^{(m)})^{j}_t=(\B_t^{(m+1)})^j.
\ee
Since for $j\leq m$ we have $(\B_t^{(m)})^j=(\B_t^{(m+1)})^j$,
\[
A \bar{N}^{(m)}_t+\B_t^{(m)}=
A \bar{N}^{(m+1)}_t+\B_t^{(m+1)}+L^{(m)}_t+Z^{(m)}_t.
\]
So, if we show that the contribution given by $Z^{(m)}$ satisfies condition {\bf (R)}, the inductive step in \eqref{indB} follows. Recall \eqref{Y}.
From \eqref{depDB}, with an abuse of notation we write
\[
(\partial_t \partial_B^{m-1} B)_j(s_{m+1} ,X_{s_{m+1}})
=
(\partial_t \partial_B^{m-1} B)_j(s_{m+1} ,X^{(j-m)}_{s_{m+1}})
\]
We recall $(H_2)$. Using the boundedness of the derivatives of $B$ we can write
\begin{multline}
\big| (\partial_t \partial_B^{m-1} B)_j(s_{m+1} ,X_{s_{m+1}}) - (\partial_t \partial_B^{m-1} B)_j(s_{m+1} ,\th_{s_{m+1}}) \big|\\ 
\leq
C |X^{(j-m)}_{s_{m+1}}-\th^{(j-m)}_{s_{m+1}}|
(1+|X^{(j-m)}_{s_{m+1}}|^{m-1}+|\th^{(j-m)}_{s_{m+1}}|^{m-1}),
\end{multline}
where $C$ depends on $T,m$ and $\kappa$.
We know from \eqref{eqmoments} that 
$\|X_{s_{m+1}}^{(j-m)}-\th_{s_{m+1}}^{(j-m)}\|_{p} \leq C_{p} {s_{m+1}}^{j-m-1/2}$. Moreover, from the Gronwall Lemma, $|\th^{(j-m)}_{s_{m+1}}|^{m-1}\leq C_\kappa$. Since $s_{m+1}\leq t\leq T$,
\[
\big\| 1+|X^{(j-m)}_{s_{m+1}}|^{m-1}+|\th^{(j-m)}_{s_{m+1}}|^{m-1} \big\|_{p} \leq C_{p,\kappa}. 
\]
We obtain
\[
\big\| (\partial_t \partial_B^{m-1} B)_j({s_{m+1}} ,X_{s_{m+1}}) - (\partial_t \partial_B^{m-1} B)_j({s_{m+1}} ,\th_{s_{m+1}}) \big\|_p\leq
C {s_{m+1}}^{j-m-1/2}
\]
Now, integrating $m+1$ times, we conclude that
\be{estZ}
\|(Z^{(m)})^j_t\|_{p} \leq C t^{m+1} t^{j-m-1/2}= C t^{j+1/2}.
\ee
We now consider the first order Malliavin derivative. Recall \eqref{depDB} and $(H_2)$. We apply the chain rule, and find that for some positive integer $m$,
 \[
|D_u (\partial_t \partial_B^{m-1} B)_j(s_{m+1} ,X_{s_{m+1}})|
\leq
C (1+| X^{(j-m)}_{s_{m+1}}|^m ) \,| D_u X^{(j-m)}_{s_{m+1}}|.
\]
Integrating in $dt$ $m+1$ times and using \eqref{eqmoments}, we obtain
\[
\|(Z^{(m)}_t)^j\|_{1,p} \leq C t^{m+1} t^{j-m-1/2}= C t^{j+1/2},
\]
so this term satisfies condition {\bf (R)} as well.
For the higher order Sobolev norms the proof is analogous. We have finally proved \eqref{indB}. For $m=n$, we have
\[
X_t-\th_t= A N_t+\B^{(n)}_t+\bar{ R}_t^{(n)}.
\]
We now prove that $\B^{(n)}_t$ satisfies condition {\bf (R)}, and this implies the statement. Recall
\[
(\B^{(n)}_t)^j=\int_0^t \dots \int_0^{s_{j -1}}  
(\partial_B^{ j-1} B)_j(s_{ j} ,X_{s_{j}})
-(\partial_B^{j-1} B)_j (s_{j},\th_{s_{j}}) d s_{ j} \dots ds_1
\]
Writing explicitly $(\partial_B^{ j-1} B)$, the Sobolev norms can now be estimated as the Sobolev norms of $(Z^{(m)}_t)^j$ in \eqref{estZ}.
\epr

\bpr\emph{(of Theorem \ref{theorem:epsilon})}
From Lemma \ref{decomposition}. For $t\in(0,1]$ we have 
\[
\chi_t  = \cal{T}_t^{-1} (X_t-\th_t) =G + \cal{T}_t^{-1} R_t,
\] 
where $R_t$ satisfies condition {\bf (R)} and $G=A \Theta$ with
\be{Theta}
\Th=\left(
\begin{array}{c}
\Th_1\\
\vdots \\
\Th_n\\
\end{array}
\right),
\quad \quad
\Th_j= \frac{N_t^j}{t^{j-1/2}} \mbox{ r.v. in } \R^d, \mbox{ for } j=1,\dots, n.
\ee
$\Th$ is a non-degenerate Gaussian r.v. in $\R^{nd}$ since its covariance can be expressed as a block matrix as
\[
Q=Cov(\Th)=\left( \frac{\Id_d}{(l+j-1)(l-1)!(j-1)!} \right)_{1\leq l,j\leq n}
\]
(cf. $(D_4)$).  We have $\|G\|_{k,p}\leq C$ for any $k\in \N,\,p\geq 2$. Since $R_t$ satisfies condition {\bf (R)},  
$\|\cal{T}_t^{-1} R_t\|_{k,p}\leq C\sqrt{t}$. This also implies $\|\chi_t \|_{k,p}\leq C$. Moreover, we have proved in Lemma \ref{nondegchain} that $\Gamma_{\chi_t}(p)\leq C$. From \eqref{deri2} with $F=\chi_t$ we get that for any multi-index $\a\in (1,\dots, dn)^k$ there exist $K>0$ constant depending on $\kappa,\lambda,T$ such that for $t\in (0,T]$
\be{ing}
|\partial^\a_z p_G(z) - \partial^\a_z p_F(z)|\leq  K  \sqrt{t} .
\ee
The variance of $G$ is $AQA^T$. We recall gradient and Hessian of the centred $nd$-dimensional Gaussian density $p_G$ with variance $\Sigma$:
\be{gh}
\nabla p_G(z)=-p_G(z)\Sigma^{-1}z,\quad\quad
\nabla^2 p_G(z)= p_G(z)(\Sigma^{-1} z z^T \Sigma^{-1} - \Sigma^{-1})
\ee
Remark, from $\chi_t= \cal{T}_t^{-1}(X_t-\th_t)$, that 
\be{cv}
\begin{split}
p_{X_t}(y)&=\frac{1}{\det \cal{T}_t }p_{\chi_t}(\cal{T}_t ^{-1}(y-\th_t)),\\
\nabla p_{X_t}(y)&=\frac{\cal{T}_t ^{-1}}{\det \cal{T}_t }  \nabla p_{\chi_t}(\cal{T}_t ^{-1}(y-\th_t)),\\
\nabla^2 p_{X_t}(y)&=\frac{\cal{T}_t ^{-1}}{\det \cal{T}_t}  \nabla^2 p_{\chi_t}\big(\cal{T}_t ^{-1}(y-\th_t)\big)\cal{T}_t^{-1}.
\end{split}
\ee
We set $\bar{y}_t=\cT_t^{-1}(y_t-\th_t)$. Using \eqref{ing}, if \eqref{condy} holds, $\bar{y}_t\rightarrow \bar{y}$ and
\[
\partial^\a_y p_{\chi_t}(\cal{T}_t ^{-1}(y_t-\th_t)) 
\rightarrow \partial^\a_z p_G(\bar{y}).
\]
Now  \eqref{cv}, \eqref{gh} with $\Sigma = AQA^T$ and \eqref{qn} give the statement.
\epr

\appendix

\section{Tools of Malliavin Calculus for density estimates}\label{app-mall}
We recall some basic notions in Malliavin calculus. Our main reference is \cite{Nualart:06}. We consider a probability space $(\Omega,\mathcal{F},\PR)$ and a Brownian motion $W=(W^1_t,...,W^d_t)_{t\geq 0}$ and the filtration $(\mathcal{F}_t)_{t \geq 0}$ generated by $W$.

For any $k \in \N$, for a multi-index $\alpha=(\alpha_1,...,\alpha_k)\in\{1,...,d\}^k$ and $(s_1,...,s_k)\in [0,T]^k$, we denote the corresponding Malliavin derivative 
$$
D^\alpha_{s_1,...,s_k} F := D^{\alpha_1}_{s_1} ... D^{\alpha_k}_{s_k} F,
$$
with $F$ random variable on $\Omega$ differentiable enough in Malliavin sense.

We denote the Malliavin-Sobolev norm of order $k\in \N $ and integrability $p\in [0,\infty)$ as
$$
\|F\|_{k,p}=[\E|F|^p+\sum_{j=1}^k \E |D^{(j)} F|^p]^\frac{1}{p}
$$
where
$$
|D^{(j)} F|=\left(\sum_{|\alpha|=j}\int_{[0,T]^j}|D^\alpha_{s_1,...,s_j}F|^2 d s_1 ... d s_j\right)^{1/2}.
$$
Clearly,
$
\|F\|_p=
\|F\|_{0,p}=\E[|F|^p]^{1/p}
$
simply stands for the $L^p$ norm as a random variable.
For $k=1$, we also use the notation $|DF|=|D^{(1)}F|$.
We denote by $\DD^{k,p}$ the space of the random variables which are $k$ times differentiable in the Malliavin sense in $L^p$, and $\DD^{k,\infty}=\bigcap_{p=1}^\infty \DD^{k,p}$. 

For a random vector $F=(F_1,...,F_n)$ in the domain of $D$, we denote its \emph{Malliavin covariance matrix} as
\[
\g_F^{i,j}= \int_{[0,T]} D_s F_i D_s F_j ds = \sum_{k=1}^d \int_0^T D^k_s F_i\times D^k_s F_j ds
\]
(in the last formula, $D^k$ denotes the derivative with respect to $W^k$. Note that this is different form $D^{(k)}$, the derivative of order $k$).
We also define
\be{gamma}
\Gamma_F(p)=1+\E\l_*(\g_F)^{-p}.
\ee
We state a known result on estimates for density and derivatives of the density of random variables. It is a version of  \cite[Theorem 2.1]{BC14}, when one takes into account also \cite[Remark 2.2 and Remark 2.4]{BC14}, with the trivial localization $\Th\equiv1$, and uses Meyer's inequality to get rid of Ornstein-Uhlenbeck operators in the bounds.
\bl{deri}
Let  $k \in \N$, and a r.v.  $F\in \DD^{k+2,p}$, with $\Gamma_F(p)<\infty$. 
Then the law of $F$ is absolutely continuous with respect to the Lebesgue measure, and for any multi-index $\alpha=(\alpha_1,...,\alpha_k)\in (1,\dots, dn)^k$,
 there exist $C, b, p \in \R^+$ such that, for every $y\in \R^n$,
\[
|\partial_x^\a p_F(y)| \leq C \Gamma_F(p)\,\|F\|_{k+2,p}\,\PR(|F|>|y|/2)^b.
\]
Moreover, if $F,G\in \DD^{k+3,p}$ are two r.v.s with $\Gamma_F(p)<\infty,\,\Gamma_F(p)<\infty$, for any multi-index $\alpha=(\alpha_1,...,\alpha_k)\in (1,\dots, dn)^k$, there exist $C, b, p \in \R^+$ such that, for every $y\in \R^n$,
\be{deri2}
|\partial_x^\a p_F(y)-\partial_x^\a p_G(y)| \leq C \Gamma_F(p)\Gamma_G(p)\,( \|F\|_{k+3,p}+\|G\|_{k+3,p})\|F-G\|_{k+2,p}.
\ee
\el

\bibliographystyle{plain} 
\begin{spacing}{0.9}
\begin{small}
\bibliography{bibliografia}
\end{small}
\end{spacing}

\end{document}